\documentclass[11pt]{article}
\usepackage{amsthm,amsmath,amssymb,anysize}
\def\qed{\hbox to 0pt{}\hfill$\rlap{$\sqcap$}\sqcup$}
\newtheorem{lemma}{Lemma}[section]
\newtheorem{theorem}[lemma]{Theorem}

\newtheorem{proposition}[lemma]{Proposition}

\newtheorem{corollary}[lemma]{Corollary}
\marginsize{3cm}{3cm}{2.3cm}{2.2cm}
\numberwithin{equation}{section}

\title{Derivations for the even part of the
 odd Hamiltonian superalgebra in positive characteristic\footnote{Supported by a NSF grant 10471096 of China and
``One Hundred Talents Program" from University of Science and Technology of China } }
 \author{Wende Liu$^{1, 2,\,}$\footnote{Corresponding author.
  Email-addresses: wendeliu@ustc.edu.cn (W. Liu), huaxiuyingnihao@163.com (X. Hua), ycsu@ustc.edu.cn (Y. Su)},
 \, Xiuying Hua$^{2}$, and Yucai Su$^{1}$\\
 \\
 $^{1}$Department of Mathematics \\
  University of Science and Technology of China, Hefei
  230026\\
  \\
   $^{2}$Department of Mathematics \\
  Harbin Normal University, Harbin
  150080
 }
\date{ }

\begin{document}
\maketitle

\begin{abstract}
\noindent In this paper we mainly study the derivations for even part of the
finite-dimensional odd Hamiltonian superalgebra $HO$ over a field of prime
characteristic. We first
give the generator set of $HO_{\overline{0}}.$ Then we determine the homogeneous derivations
from $HO_{\overline{0}}$ into
 $W_{\overline{0}},$ the even part of the generalized Witt superalgebra $W$. Finally, we
 determine the derivation algebra and outer derivation algebra of
 $HO_{\overline{0}}.$\\


  \end{abstract}
\vspace{1cm}

\noindent \textbf{0. Introduction}\vspace{0.5cm}

\noindent The theory of Lie superalgebras has seen a significant
development (c.f. \cite{sc}). For example, V. G. Kac has completed the
classification of finite-dimensional simple Lie superalebras and the classification of
infinite-dimensional simple linearly compact Lie superalgebras over algebraically
closed
fields of characteristics zero (see \cite{k1,k2}).
But there are not so plentiful results
on modular Lie superalgebras (that is, Lie superalgebras
 over fields of finite characteristic). For example, the classification problem is still
 open for finite-dimensional simple Lie superalgebras. As far as we know, the
 references \cite{kl,p} should be the earliest
 papers on modular Lie superalgebras.

In this paper we shall discuss the derivations  for  the even part of
the finite-dimensional odd Hamiltonian superalgebra $HO$ over a field of prime characteristic.
Our work is essentially motivated by the work on modular Lie algebras of Cartan type
(see \cite{c,st,sf}).  In \cite{mz,wz,zz} the superderivation algebras were
determined for the finite-dimensional modular Lie superalgebras of
 Cartan type
$W,$ $S,$ $H,$ and $K.$ In \cite{lzw} the superderivation
algebra was determined for the finite-dimensional odd Hamiltonian superalgebra
$HO.$ Note that the
derivations of the even
parts have  been studied  for the Lie superalgebras of Cartan type $W,$  $S,$
and $H$ (see \cite{lsz,lz1}). However, in contrast to the setting of $H$,
 we determine completely the derivation space from the even part of $HO$ into the even part of  $W,$
and  the derivation algebra of the even part of $HO$, also.

This paper is arranged as follows. In Section 1 we introduce the
necessary notations, definitions and known results. In Section 2
we first give the generating set of the even part of the odd
Hamiltonian superalgebra. Then  we determine the nonnegative
$\mathbb{Z}$-degree derivations from the even part of the odd
Hamiltonian superalgebra into the even part of the generalized
Witt superalgebra. In Section 3 we determine the  negative
$\mathbb{Z}$-degree derivations from the even part of the odd
Hamiltonian superalgebra into the even part of the generalized
Witt superalgebra. In Section 4 we determine completely the
derivation algebra and outer derivation algebra of the even part
of the odd Hamiltonian superalgebra.

\section{Preliminaries}

Let $\mathbb{Z}_{2}=\{\overline{0}, \overline{1}\}$ be the field
of two elements. For a   vector superspace
$V=V_{\overline{0}}\oplus V_{\overline{1}}, $ we denote by
$\mathrm{p}(a)=\theta$ the \textit{parity of a homogeneous
element} $a\in V_{\theta}, \theta\in \mathbb{Z}_{2}.$
 We assume throughout  that the notation   $\mathrm{p}(x)$ implies
that
  $x$ is a $\mathbb{Z}_2$-homogeneous element.

Let $\frak{g}$ be a  Lie algebra and $V$   a $\frak{g}$-module. A
linear mapping $D:\frak{g}\rightarrow V$ is called a
\textit{derivation } from $\frak{g}$  into $V$ if $D(xy)=x\cdot
D(y)-y\cdot D(x)$ for all $x,y\in \frak{g}.$ A derivation
$D:\frak{g}\rightarrow V$ is called \textit{inner} if there is
$v\in V $ such that $D(x)=x\cdot v$ for all $x\in \frak{g}.$
Following \cite[p. 13]{sf}, denote by $\mathrm{Der}(\frak{g},V)$
the \textit{derivation space} from $\frak{g}$  into $V.$  Then
$\mathrm{Der}(\frak{g},V)$ is a $\frak{g}$-submodule of
$\mathrm{Hom}_{\mathbb{F}}(\frak{g},V).$ Assume in addition that
$\frak{g}$ and $V$ are finite-dimensional and that
$\frak{g}=\oplus_{r\in \mathbb{Z}} \frak{g}_{r}$ is
$\mathbb{Z}$-graded and $V=\oplus_{r\in\mathbb{ Z}} V_{r}$ is a
$\mathbb{Z}$-graded $\frak{g}$-module. Then
$\mathrm{Der}(\frak{g},V)=\oplus_{r\in \mathbb{Z}}
\mathrm{Der}_{r}(\frak{g},V)$ is a $\mathbb{Z}$-graded
$\frak{g}$-module  by setting
$$\mathrm{Der}_{r}(\frak{g},V):=\{D\in
\mathrm{Der}(\frak{g},V)\mid D(\frak{g}_{i})\subset V _{r+i} \
\mbox{for all} \ i \in\mathbb{Z}\}.$$ In the case of $V=\frak{g},$
the \textit{derivation algebra} $\mathrm{Der}(\frak{g})$ coincides
with $\mathrm{Der}(\frak{g},\frak{g}) $ and
$\mathrm{Der}(\frak{g}) =\oplus _{r\in \mathbb{Z}}
 \mathrm {Der}_{r}(\frak{g})$ is a $\mathbb{Z}$-graded Lie algebra.
If $\frak g=\oplus_{-r\leq i\leq s}\frak g_i$ is a $\mathbb
    Z$-graded Lie  algebra, then $\oplus_{-r\leq i\leq 0}\frak
    g_i$ is called the \textit{top of} $\frak g$ (with respect to the
    gradation).

 In the sequel we assume throughout that the underlying filed
  $\mathbb{F}$ is  of characteristic $p>3.$
  Let us  review the notions of  modular Lie
superalgebras of Cartan-type $W$ and $S$ and their gradation
structures.
 In addition to the standard notation
$\mathbb {Z},$ we write  $\mathbb {N}$ for the set of positive
integers, and ${\mathbb {N}}_0$ for the set of nonnegative
integers. Henceforth, we will let  $n$  denote fixed integers in
$\mathbb {N}\setminus \{1,2\} $ without notice. For $\alpha = (
\alpha _1,\ldots,\alpha _n ) \in \mathbb {N}_0^n,$ we put
$|\alpha| =\sum_{i=1}^n\alpha _i.$ Let $\mathcal{O}(n)$ denote the
\textit{divided power algebra over} $\mathbb{F}$ with an ${\mathbb
F}$-basis $\{ x^{( \alpha ) }\mid \alpha \in \mathbb{N}_0^n \}. $
For $\varepsilon _i=( \delta_{i1}, \ldots,\delta _{in}),$ we
abbreviate $x^{( \varepsilon _i)}$ to $x_i,$ $i=1,\ldots,n.$ Let
$\Lambda (n)$ be the \textit{exterior superalgebra over}
$\mathbb{F}$ in $n$ variables $x_{n+1},\ldots ,x_{2n}.$ Denote
the tensor product by $\mathcal{O}(n,n)
=\mathcal{O}(n)\otimes_{\mathbb{F}} \Lambda(n).$ Obviously,
$\mathcal{O}(n,n)$ is an associative superalgebra with a
$\mathbb{Z}_2$-gradation induced by the trivial
$\mathbb{Z}_2$-gradation of $\mathcal{O}(n)$ and the natural
$\mathbb{Z}_2$-gradation of $\Lambda (n).$ Moreover,
$\mathcal{O}(n,n)$ is super-commutative. For $g\in
\mathcal{O}(n),f\in \Lambda(n),$ we write $gf $  for $ g\otimes
f.$ The following formulas hold in $\mathcal{O}(n,n):$
$$
x^{(\alpha) }x^{(\beta) }=\binom{\alpha +\beta }{ \alpha}
 x^{( \alpha +\beta)}\quad\mbox{for}\  \alpha,\beta \in
{\mathbb N}_0^n;
$$
$$
x_kx_l=-x_lx_k\quad\mbox{for}\ k,l=n+1,\ldots,2n;
$$
$$
x^{( \alpha ) }x_k=x_kx^{( \alpha ) }\quad\mbox{for}\ \alpha \in
\mathbb{N}_0^n, k=n+1,\ldots,2n,
$$
where $\binom{ \alpha +\beta} {\alpha}:=\prod_{i=1}^n\binom{
\alpha _i+\beta _i}{ \alpha _i}.$ Put $Y_0:=\{ 1,2,\ldots,n \},$
$Y_1:=\left\{n+1,\ldots,2n\right\} $ and $Y:=Y_0\cup Y_1.$
Let
$$\mathbb{B}_k:=\left\{ \langle i_1,i_2,\ldots
,i_k\rangle\mid n+1\leq i_1<i_2<\cdots <i_k\leq 2n\right\}  $$
be the set of $k$-tuples of strictly increasing integers between $n+1$
and $2n,$ and put $
\mathbb{B}:=\mathbb{B}(n):=\bigcup\limits_{k=0}^n\mathbb{B}_k,$
where $\mathbb{B}_0:=\emptyset.$ Put $\mathbb{B}^0:=\{u\in
\mathbb{B} \mid  |u| \, \mbox{even}\}$ and
 $\mathbb{B}^1:=\{u\in \mathbb{B} \mid  |u| \,\mbox{odd}\},$ where for $u=\langle i_1,i_2,\ldots
,i_k\rangle \in \mathbb{B}_k,$
 $ |u| :=k,| \emptyset | :=0$,
$x^\emptyset :=1.$ For $u=\langle i_1,i_2,\ldots
,i_k\rangle \in \mathbb{B}_k,$ we set $x^u:=x_{i_1}x_{i_2}\cdots x_{i_k};$  we also use   $u$ to stand for the set $\{i_1,i_2,\ldots
,i_k\} $ if no confusion occurs.
 Clearly, $\left\{ x^{\left( \alpha \right) }x^u\mid\alpha
\in \mathbb{N}_0^n,u\in \mathbb{B}  \right\} $ constitutes an
$\mathbb{F}$-basis of $\mathcal{O} \left( n,n\right). $ Let
$\partial_1,\partial_2,\ldots,\partial_{2n}$ be the linear transformations of
$\mathcal{O} \left( n,n\right) $ such that
\[
\partial_r ( x^{\left( \alpha \right) }x^u ) =\left\{
\begin{array}{l}
x^{\left( \alpha -\varepsilon _r\right) }x^u,\quad \quad \quad
\quad r\in Y_0
\\
x^{\left( \alpha \right) }\cdot \partial x^u/\partial x_r,\quad
\quad r\in Y_{1.}
\end{array}
\right.
\]
Then $\partial_1,\partial_2,\ldots,\partial_{2n}$ are superderivations of the
superalgebra $\mathcal{O} \left( n,n\right).$ Let
\[
W\left(n,n\right) =\Big\{ \sum\limits_{r\in Y}f_r \partial_r \mid f_r\in
\mathcal{O} \left( n,n\right),r\in Y\Big\}.
\]
Then $W\left( n,n\right) $ is a Lie superalgebra, which is
contained in ${\rm Der} ( \mathcal{O} \left(n,n\right)). $
Obviously, ${\rm p}( \partial_i) =\mu (i),$ where
\[
\mu \left( i\right) :=\left\{
\begin{array}{l}
\overline{0},\quad \quad i\in Y_0 \\
\overline{1},\quad \quad i\in Y_{1.}
\end{array}
\right.
\]
One may verify that
\[[fD, gE]=fD(g)E-(-1)^{{\rm
p}(fD){\rm p}(gE)}gE(f)D+(-1)^{{\rm p}(D){\rm p}(g)}fg[D, E]\] for
$ f,g\in \mathcal{O} (n,n),$ $ D,E\in {\rm Der}( \mathcal{O} (n,n
)).$  Let
\[
\underline{t}:=\left( t_{1},t_{2},\ldots,t_n\right) \in
\mathbb{N}^n,\quad \pi:=\left( \pi _1,\pi _2,\ldots,\pi _n\right)
\]
where $\pi _i:=p^{t_i}-1,i\in Y_0.$ Let
$\mathbb{A}:=\mathbb{A}\left( n;\underline{t}\right) =\left\{
\alpha \in \mathbb{N}_0^n\mid\alpha _i\leq \pi _i,i\in Y_0
\right\}.$ Then
\[
\mathcal{O} \left( n,n;\underline{t}\right) :={\rm
span}_{\mathbb{F}}\left\{ x^{\left( \alpha \right) }x^u \mid
\alpha \in \mathbb{A},u\in \mathbb{B} \right\}
\]
is a finite-dimensional subalgebra of $\mathcal{O}
\left(n,n\right) $ with   a natural ${\mathbb Z}$-gradation
$\mathcal{O} \left(n,n;\underline{t}\right)=\bigoplus _{r=0}^{\xi}
\mathcal{O}(n,n;\underline{t})_{r}$ by putting
        $$\mathcal{O}(n,n;\underline{t})_{r}:=
        {\rm span} _{\mathbb F}\{ x^{(\alpha)} x^{u}\mid |\alpha|+|u|=r \},
        \quad \xi:=|\pi|+n.$$
Set
\[
W\left( n,n;\underline{t}\right):=\Big\{ \sum\limits_{r\in Y}
f_r\partial_r \mid f_r\in \mathcal{O} \left(
n,n;\underline{t}\right),r\in Y\Big\}.
\]
Then $W\left( n,n;\underline{t}\right) $ is  a finite-dimensional
simple Lie superalgebra (see \cite{z}). Obviously, $W (
n,n;\underline{t} ) $ is a free $ \mathcal{O} \left(
n,n;\underline{t}\right)$-module with $ \mathcal{O} \left(
n,n;\underline{t}\right) $-basis $ \{ \partial_r\mid r\in Y \}.$ We note
that $W ( n,n;\underline{t} ) $ possesses a \textit{standard}
$\mathbb F$-\textit{basis} $\{x^{(\alpha)}x^{u}\partial_r\mid \alpha\in
\mathbb{A}, u\in \mathbb{B}, r\in Y\}.$  The
$\mathbb{Z}$-gradation of    $\mathcal{O}(n,n;\underline{t})$
induces naturally   a $\mathbb{Z}$-gradation structure   of
  $W(n,n;\underline{t})
=\oplus_{i=-1}^{\xi-1} W(n,n;\underline{t})_{i},$
  where
$$W(n,n;\underline{t})_{i} :={\rm span}_{\mathbb{F}}\{f\partial_{s}\mid
s\in Y,\ f\in \mathcal{O}(n,n;\underline{t})_{i+1}\}.$$ Put
\[
i^{\prime }=\left\{
\begin{array}{l}
i+m,\quad \quad i\in Y_0 \\
i-m,\quad \quad i\in Y_{1.}
\end{array}
\right.
\]
 Define
  $\mathrm{T_H}:\mathcal{O}(n,n;\underline{t})\rightarrow
 W(n,n;\underline{t})$ such that
 \[
 \mathrm{T_H}(a):=\sum_{i=1}^{2n}(-1)^{\mu(i)\mathrm{p}(a)}\partial_i(a)\partial_{i'}
  \quad\mbox{for all}\  a\in
 \mathcal{O}(n,n;\underline{t}).
 \] Then $\mathrm{T_H}$ is an odd linear mapping. For $a,\,b\in
\mathcal{O}(n,n;\underline{t}),$ one may easily verify that
\begin{equation}\label{he1.1}
[\mathrm{T_H}(a),\mathrm{T_H}(b)]=\mathrm{T_H}(\mathrm{T_H}(a)(b))
\end{equation}
(see \cite[Proposition 1, Lemma 2]{lzw}).
Put \[
HO(n,n;\underline{t}):=\{\mathrm{T_H}(a)\,|\,a\in\mathcal{O}(n,n;\underline{t})\}.
\]
Then (\ref{he1.1}) shows that $HO(n,n;\underline{t})$ is a finite-dimensional Lie superalgebra. Following \cite{k2}, we call this Lie
superalgebra
  the
odd Hamiltonian superalgebra (see also \cite{lzw}).  Put
\[
\overline{HO}(n,n;\underline{t}):=\overline{HO}(n,n;\underline{t})_{\overline{0}}
+\overline{HO}(n,n;\underline{t})_{\overline{1}},\] where
\[\overline{HO}(n,n;\underline{t})_\alpha
:=\Big\{\sum_{i=1}^{2n}a_i\partial_i\in
W(n,n;\underline{t})_\alpha\,\Big|\,
\partial_i(a_{j'})=(-1)^{(i,j,\alpha)}\partial_j(a_{i'}),\,\, i,j\in Y\Big\}  \]
where $\alpha \in \mathbb{Z}_2,$
$(i,j,\alpha):=\mu(i)\mu(j)+(\mu(i)+\mu(j))(\alpha+\overline{1}).$
Note that $HO(n,n;\underline{t})$ and
$\overline{HO}(n,n;\underline{t})$ are all $\mathbb{Z}$-graded
subalgebras of $W(n,n;\underline{t}).$

In the sequel the even parts of $W(n,n;\underline{t}),$ $HO(n,n;\underline{t})$ and
$\overline{HO}(n,n;\underline{t})$ will be denoted by
$\mathcal{W}$, $\mathcal{HO}$ and $\overline{\mathcal{HO}},$ respectively.

Let $\mathcal{G}:=\mathrm{span}_\mathbb{F}\{x^u\partial_r\,|\,r\in Y,u\in
\mathbb{B},\mathrm{p}(x^u\partial_r)=\overline{0}\}.$
Then $C_{\mathcal{W} } (\mathcal{HO} _{-1})=\mathcal{G};$ in particular,
$C_{\mathcal{HO} } (\mathcal{HO} _{-1})\subset\mathcal{G}.$ Note that $\mathcal{G}$
is a $\mathbb{Z}$-graded subalgebra of
$\mathcal{W}.$ We need the following lemmas.

\begin{lemma}\label{ht1.2}
 Let $\phi\in
 \mathrm{Der}(\mathcal{HO},\mathcal{W})$ satisfy  $\phi(\mathcal{HO}_{-1})=0.$ Let $E\in
 \mathcal{HO}.$ Then $[E,\mathcal{W}_{-1}]\subseteq\ker
 \phi$ if and only if $\phi(E)\in \mathcal{G}.$
\end{lemma}
\begin{proof}
 Note that $\mathcal{HO}$ is a $\mathbb{Z}$-graded subalgebra of
 $\mathcal{W}$ and
 $\mathcal{HO}_{-1}=\mathcal{W}_{-1}.$
  This is a direct consequence of  \cite[Lemma 2.1.1]{lz1}.
\end{proof}

\begin{lemma}\label{ht1.3} (\cite[Lemma 2.1.3]{lz1})
Let $V$ be a vector space over $\mathbb{F}$ and
$v_1,v_2,\dots,v_k\in V.$ Let $A_i\in \mathrm{End}_\mathbb{F}V.$
Suppose   there are  $B_i\in \mathrm{End}_\mathbb{F}V $
such that

\begin{enumerate}

\item [$\mathrm{(i)}$]  $A_iB_iA_i=A_i $ for $1\leq i\leq k;$

\item [$\mathrm{(ii)}$] $A_1,\,A_2,\,\dots,\,A_k $ commute mutually;

\item [$\mathrm{(iii)}$] $A_i(v_j)=A_j(v_i)$ for $ 1\leq i,j\leq k;$

\item [$\mathrm{(iv)}$] $A_iB_i(v_i)=v_i,$ $ A_iB_j=B_jA_i $ for $1\leq i,j\leq
k,$ $i\neq j$.
\end{enumerate}
 Then there is  $v\in V $ such that $A_i(v)=v_i $ for $1\leq i\leq k.$
\end{lemma}

\begin{lemma}\label{ht1.4} Let $\phi \in \mathrm{Der}_t(\mathcal{HO},\mathcal{W}) $ where
$t\geq 0.$ Then there is $E \in\mathcal{W}_t $
such that $(\phi-\mathrm{ad}E)(\mathcal{HO}_{-1})=0.$
\end{lemma}
\begin{proof}
 Apply \cite[Proposition 2.1.6]{lz1}.\end{proof}

\section{Nonnegative $\mathbb{Z}$-degree derivations   from
$\mathcal{HO}$ into $\mathcal{W}$}

In order to compute the derivations from $\mathcal{HO}$ into
$\mathcal{W}$, we first give the generating set of $\mathcal{HO},$ which will be frequently used in the sequel.
Let $$N:=\{\mathrm{T_H}(x_kx_lx_q)\,|\,k,l,q\in Y_1\},$$
$$M:=\{\mathrm{T_H}(x^{(q_i\varepsilon_i)}x_k)\,|\,i\in Y_0,0\leq
q_i\leq \pi_i,k\in Y_1\}.$$

\begin{proposition}\label{ht2.1}
  $\mathcal{HO}$ is generated by $M\bigcup N.$
\end{proposition}

\begin{proof} Let $L$  be the subalgebra of $\mathcal{HO}$ generated  by $M\cup
N.$
We proceed by induction on $|\alpha|+|u|$ to show that
$\mathrm{T_H}(x^{(\alpha)}x^u)\in L$ for all $\alpha \in \mathbb{A}$
and $u\in \mathbb{B}^1.$ When $|\alpha|+|u|=1$, we have $\alpha=0$ and therefore,
$\mathrm{T_H}(x_k)\in M\subset L$ for $k\in Y_1.$ When $|\alpha|+|u|=2$,
$\mathrm{T_H}(x_ix_k)\in M\subset L$ for $i\in Y_0,\,k\in Y_1.$ When
$|\alpha|+|u|=3$, we have
$$ [\mathrm{T_H}(x_ix_{j'}),\,\mathrm{T_H}(x^{(2\varepsilon_j)}x_l)]
 =-(1+\delta_{ij})\mathrm{T_H}(x^{(\varepsilon_i+\varepsilon_j)}x_l)
 \in L
 $$ for $ i,j\in Y_0 $ and $ l\in Y_1 $ with $l\neq
 i';$
 also,
$$[\mathrm{T_H}(x_jx_{i'}),\,\mathrm{T_H}(x^{(2\varepsilon_
i)}x_{i'})]=-\mathrm{T_H}(x^{(\varepsilon_
i+\varepsilon_j)}x_{i'})\in L \quad\mbox{for}\, i,j\in Y_0.$$
Summarizing, $\mathrm{T_H}(x^{(\varepsilon_i+\varepsilon_j)}x_k)\in L$ for all
$k\in Y_1,$ $i,j\in
 Y_0.$
  Since $\mathrm{T_H}(x_kx_lx_q)\in N\subset L$ for $k,l,q\in Y_1,$
  our assertion holds for $|\alpha|+|u|=3.$
Now Suppose $|\alpha|+|u|>3.$
Direct computation shows that
$$[\mathrm{T_H}(x^{(\pi_1\varepsilon_1)}x_{2'}),
\,\mathrm{T_H}(x^{(\pi_2\varepsilon_2)}x_k)]=
-\mathrm{T_H}(x^{(\pi_1\varepsilon_1+(\pi_2-1)\varepsilon_2)}x_k)\in
L \quad \mbox{for all}\, k\in Y_1\setminus 1';$$
$$[\mathrm{T_H}(x^{(\pi_1\varepsilon_1+(\pi_2-1)\varepsilon_2)}x_k),\,\mathrm{T_H}(x^{(\varepsilon_
{k'}+\varepsilon_2)}x_k)]=(1+\delta_{k,2'})
\mathrm{T_H}(x^{(\pi_1\varepsilon_1+\pi_2\varepsilon_2)}x_k)\in
L\quad\mbox{for all}\, k\in Y_1\setminus 1';$$
$$[\mathrm{T_H}(x^{(2\varepsilon_1)}x_{1'}),
\,\mathrm{T_H}(x^{(\pi_2\varepsilon_2)}x_{1'})]=
\mathrm{T_H}(x^{(\varepsilon_1+\pi_2\varepsilon_2)}x_{1'})\in L;$$
$$[\mathrm{T_H}(x^{(\pi_1\varepsilon_1)}x_{1'}),\,\mathrm{T_H}(x^{(\varepsilon_1+\pi_2\varepsilon_2)}x_{1'})]=
-2\mathrm{T_H}(x^{(\pi_1\varepsilon_1+\pi_2\varepsilon_2)}x_{1'})\in
L.$$
Since $\mathrm{char}\mathbb{F} >3,$
  $\mathrm{T_H}(x^{(\pi_1\varepsilon_1+\pi_2\varepsilon_2)}x_k)\in
L$ for $ k\in Y_1.$  Similarly, one may  show inductively that
   $\mathrm{T_H}(x^{(\pi)}x_k)\in L$ for all $ k\in
Y_1.$ Clearly,
$$[\mathrm{T_H}(x^{(\pi)}x_k),\,\mathrm{T_H}(x_{i'})]
=\mathrm{T_H}(x^{(\pi-\varepsilon_i)}x_k)\in
L\quad \mbox{for all}\, i\in Y_0,\,k\in Y_1.$$
 Therefore,
 $\mathrm{T_H}(x^{(\alpha)}x_k)\in L $ for all $
 k\in Y_1,\,\alpha\in \mathbb{A}.$
 Thus one may assume that  $|u|\geq3.$
Let $x^u=x_{i_{1'}}x_{i_{2'}}\cdots x_{i_{s'}},$ where $i_{j'}\in
 Y_1 $ and $s\geq 3$ is odd.
If there is some $j'\in \{1',2',\dots,s'\} $ such that  $\alpha_{i_j}<\pi_{i_j},$
without loss of generality, assume that  $j'\in \{1',2',\dots,{(s-2)}'\}.$
By inductive hypothesis,
$$[\mathrm{T_H}(x^{(\alpha+\varepsilon_{i_j})}x_{i_{1'}}\cdots
x_{i_{(s-2)'}}),\,\mathrm{T_H}(x_{i_{j'}}x_{i_{(s-1)'}}x_{i_{s'}})]
=\mathrm{T_H}(x^{(\alpha)}x^u)\in L;$$
if $\alpha_{i_{r}}=\pi_{i_{r}}$ for all $1 \leq r\leq s $ and  there is some $j'\in
\{{(s+1)}',\dots,n'\} $ such that  $\alpha_{i_j}<\pi_{i_j},$
then by inductive hypothesis,
$$[\mathrm{T_H}(x^{(\alpha+\varepsilon_{i_j})}x_{i_{1'}}\cdots
x_{i_{(s-2)'}}),\,\mathrm{T_H}(x_{i_{j'}})]=\mathrm{T_H}(x^{(\alpha)}x_{i_{1'}}\cdots
x_{i_{(s-2)'}})\in L.$$
In the first case, we have shown that
$\mathrm{T_H}(x^{(\varepsilon_{i_1})}x_{i_{1'}}x_{i_{(s-1)'}}x_{i_{s'}})\in
L.$ Thus, in the second case we have
$$[\mathrm{T_H}(x^{(\alpha)}x_{i_{1'}}\cdots
x_{i_{(s-2)'}}),\,\mathrm{T_H}(x^{(\varepsilon_{i_1})}x_{i_{1'}}x_{i_{(s-1)'}}x_{i_{s'}})]
=-2\mathrm{T_H}(x^{(\alpha)}x^u)\in L.$$
Now it suffices to show that
$\mathrm{T_H}(x^{(\pi)}x^u)\in L $ for $u\in
\mathbb{B}^1.$
 Use induction on $|u|$. For $|u|=1$ we have proved the conclusion.
For  $|u|=3$,
$$[\mathrm{T_H}(x^{(\pi)}x_k),\,\mathrm{T_H}(x^{(\varepsilon_{k'})}x_kx_rx_q)]
 =-2\mathrm{T_H}(x^{(\pi)}x_kx_rx_q)\in L\quad \mbox{for all}\, k ,r,q\in Y_1.$$
 Suppose $|u|>3 $ and write $u=v\dotplus\omega,$
 where $v=<k,l,r>.$ By the argument above,
 $\mathrm{T_H}(x^{(\varepsilon_{r'})}x_rx^\omega)\in
 L.$ Consequently,
$$[\mathrm{T_H}(x^{(\pi)}x_kx_lx_r),\,\mathrm{T_H}(x^{(\varepsilon_{r'})}x_rx^\omega)]
=-2\mathrm{T_H}(x^{(\pi)}x^u)\in L.$$ The proof is complete.\end{proof}

In order to determine the nonnegative $\mathbb{Z}$-degree derivations from $\mathcal{HO}$
into $\mathcal{W},$ we shall show that if such two derivations coincide on the top of
$\mathcal{HO}$ then the difference vanishes. Of course we shall consider the generators of
$\mathcal{HO}.$ We first consider the subset $M$.

For convenience, we put  $\Delta_{i}:=x_{i'}\partial_{i'}-x_i\partial_i= \mathrm{T_H}(x_ix_{i'})$ for $i\in
Y_0.$

\begin{lemma} \label{ht3.1} Let $\phi\in
\mathrm{Der}_t(\mathcal{HO},\mathcal{W})$ with $t\geq 0.$
If $\phi({\mathcal{HO}}_{-1}\bigoplus{\mathcal{HO}}_0)=0,$
then $$\phi(\mathrm{T_H}(x^{(a\varepsilon_i)}x_k))=0
\quad \mbox{for}\, \,  i\in
Y_0,\,k\in Y_1,\,0\leq a\leq \pi_i.$$
\end{lemma}
\begin{proof}
For $a=2$, by Lemma \ref{ht1.2} one may assume that
$$\phi(\mathrm{T_H}(x^{(2\varepsilon_i)}x_{i'}))=\sum_{r\in
Y,\,u\in \mathbb{B}}c_{ur}x^u\partial_r\quad \mbox{where}\; c_{ur}\in
\mathbb{F}.$$ Direct computation shows that
$[\mathrm{T_H}(x^{(2\varepsilon_i)}x_{i'}),\,\Delta_{i}]
=\mathrm{T_H}(x^{(2\varepsilon_i)}x_{i'}).$
Applying $\phi$ yields that
\[ \Big[\sum_{r\in Y,\,u\in
\mathbb{B}}c_{ur}x^u\partial_r,\,\Delta_{i}\Big]=\sum_{r\in
Y,\,u\in \mathbb{B}}c_{ur}x^u\partial_r.\]
Consequently,
\[\sum_{u\in
\mathbb{B}}(c_{ui'}x^u\partial_{i'}-c_{ui}x^u\partial_i)-\delta_{i'\in
u}\sum_{r\in Y,\,u\in \mathbb{B}}c_{ur}x^u\partial_r=\sum_{r\in Y,\,u\in
\mathbb{B}}c_{ur}x^u\partial_r,\]
where $\delta_{i'\in u}$ means it is $1$ if $i'\in u$ or $0$ otherwise.
A comparison of coefficients shows that
\begin{equation}
(1+\delta_{i'\in
u})c_{ur}=0\quad\mbox{for}\,r\in Y\setminus\{i,i'\},\, u\in \mathbb{B};\label{he3.1}
\end{equation}
\begin{equation}
(2+\delta_{i'\in u})c_{ui}=0\quad\mbox{for}\, u\in \mathbb{B};\label{he3.2}
\end{equation}
 \begin{equation}
 \delta_{i'\in
u}c_{ui'}=0\quad\mbox{for}\, u\in \mathbb{B}.\label{he3.3}
\end{equation}
Noticing that $\mathrm{char}\mathbb{F}>3,$ by
 (\ref{he3.1}) and (\ref{he3.2}) one gets $c_{ur}=0$ for $ r\in Y\setminus i'.$
Thus,
\[\phi(\mathrm{T_H}(x^{(2\varepsilon_i)}x_{i'}))=\sum_{u\in
\mathbb{B}^1}c_{ui'}x^u\partial_{i'}.\]
For any fixed $u\in \mathbb{B}^1,$ if $i'\in u,$
then (\ref{he3.3}) yields $c_{ui'}=0;$ if $i'\not\in u,$ since $|u|\geq
2$, there is $l\in u$ such that
 $[\mathrm{T_H}(x^{(2\varepsilon_i)}x_{i'}),\,\Delta_{l'}]=0.$
Applying $\phi$, one may obtain that
$-\sum_{i'\not\in u\in \mathbb{B}^1}c_{ui'}x^u\partial_{i'}=0.$
It follows that $c_{ui'}=0.$ Consequently,
$\phi(\mathrm{T_H}(x^{(2\varepsilon_i)}x_{i'}))=0.$

In the following we use induction on $a\geq 3$ to show that
$\phi(\mathrm{T_H}(x^{(a\varepsilon_i)}x_{i'}))=0.$  Assume our assertion holds
for  $a-1.$ By Lemma \ref{ht1.2}  one may assume that
 \[\phi(\mathrm{T_H}(x^{(a\varepsilon_i)}x_{i'}))=\sum_{r\in
Y,\,u\in \mathbb{B}}c_{ur}x^u\partial_r\quad \mbox{where}\; c_{ur}\in
\mathbb{F}.\] Note that
$$[\mathrm{T_H}(x^{(a\varepsilon_i)}x_{i'}),
\,\Delta_{i}]=(a-1)\mathrm{T_H}(x^{(a\varepsilon_i)}x_{i'}).$$
Applying $\phi,$ one gets
\[
\sum_{u\in \mathbb{B}}(c_{ui'}x^u\partial_{i'}-c_{ui}x^u\partial_i)-\delta_{i'\in
u}\sum_{r\in Y,\,u\in \mathbb{B}}c_{ur}x^u\partial_r=(a-1)\sum_{r\in Y,\;
u\in \mathbb{B}}c_{ur}x^u\partial_r.
\]
A comparison of coefficients yields
\begin{equation}
(a-1+\delta_{i'\in u})c_{ur}=0 \quad\mbox{for}\,\, r\in Y\setminus\{i,i'\},\, u\in \mathbb{B};\label{he3.4}
\end{equation}
\begin{equation}
(a+\delta_{i'\in u})c_{ui}=0 \quad\mbox{for}\,\;  u\in \mathbb{B};  \label{he3.5}
\end{equation}
\begin{equation}
(a-2+\delta_{i'\in u})c_{ui'}=0 \quad\mbox{for}\,\;   u\in \mathbb{B}.
\label{he3.6}
\end{equation}
The following discuss is divided into two parts.\\

 \noindent Part I: $a\equiv 0\pmod{p}.$
By (\ref{he3.6}), $c_{ui'}=0 $ and therefore, \[\phi
(\mathrm{T_H}(x^{(a\varepsilon_i)}x_{i'}))=\sum_{i'\in u\in
\mathbb{B}\atop r\in Y\setminus i'}c_{ur}x^u\partial_r+\sum_{i'\not\in v\in
\mathbb{B}\atop r\in Y\setminus i'}c_{vr}x^v\partial_r.\]
 By (\ref{he3.4}) and (\ref{he3.5}),
$$\phi(\mathrm{T_H}(x^{(a\varepsilon_i)}x_{i'}))=\sum_{i'\in u\in\mathbb{B}\atop
r\in Y\setminus \{i,i'\}}c_{ur}x^u\partial_r+\sum_{i'\not\in v\in
\mathbb{B}^{0}}c_{vi}x^v\partial_i.$$
For any fixed $u $ with $i'\in u,$
since $|u|\geq 3,$ there is $l\in u$ such  that $l\neq i'.$
Compute
$[\mathrm{T_H}(x^{(a\varepsilon_i)}x_{i'}),\,\Delta_{l'})]=0.$
Applying $\phi$ to this equation, one gets
$$\sum_{i'\in u\in
\mathbb{B}}(c_{ul}x^u\partial_l-c_{ul'}x^u\partial_{l'})-\sum_{i'\in u\in
\mathbb{B}\atop r\in Y\setminus \{i,i'\}}c_{ur}x^u\partial_r-\delta_{l\in
v}\sum_{i'\not\in v\in \mathbb{B}^0}c_{vi}x^v\partial_i=0,
$$
and hence $c_{ur}=0$ for $r\in Y\setminus l.$ Thus
\[\phi
(\mathrm{T_H}(x^{(a\varepsilon_i)}x_{i'}))=\sum_{i'\in u\in
\mathbb{B}^1}c_{ul}x^u\partial_l+\sum_{i'\not\in v\in
\mathbb{B}^0}c_{vi}x^v\partial_i.
\]
Given $u $ satisfying $i',l\in u,$
since $|u|\geq 3,$ there is $q\in u$ such that $q\neq i',l.$
Note that
$[\mathrm{T_H}(x^{(a\varepsilon_i)}x_{i'}),\,\Delta_{q'}]=0.$
Applying $\phi,$ one computes
$$-\sum_{i'\in u\in \mathbb{B}^1}c_{ul}x^u\partial_l-\delta_{q\in
v}\sum_{i'\not\in v\in \mathbb{B}^0}c_{vi}x^v\partial_i=0.
$$
Hence $c_{ul}=0.$
Given $v$ satisfying $i'\not\in v,$ by the equation above, one gets
$$\phi
(\mathrm{T_H}(x^{(a\varepsilon_i)}x_{i'}))=\sum_{q,i'\not\in v\in
\mathbb{B}^0}c_{vi}x^v\partial_i.
$$
Since $|v|\geq 3,$ there is $r\in v$ such that
 $
 [\mathrm{T_H}(x^{(a\varepsilon_i)}x_{i'}),\,\Delta_{r'}]=0.
 $
Applying $\phi$, we have
$$-\sum_{q,i'\not\in v\in \mathbb{B}^{0}}c_{vi}x^v\partial_i=0.$$
This implies that $c_{vi}=0.$
As a result,
$\phi(\mathrm{T_H}(x^{(a\varepsilon_ i)}x_{i'}))=0.$\\

\noindent Part II: $a\not\equiv 0\pmod{p}.$ Here we proceed in four steps.\\

\noindent Step (i):
$a\equiv 1\pmod{p}.$ By (\ref{he3.5}),  $c_{ui}=0.$
Then by (\ref{he3.4}) and  (\ref{he3.6}),
$$\phi(\mathrm{T_H}(x^{(a\varepsilon_i)}x_{i'}))=\sum_{i'\in u\in
\mathbb{B}^1}c_{ui'}x^u\partial_{i'}+\sum_{i'\not\in v\in \mathbb{B}\atop
r\in Y\setminus \{i,i'\}}c_{vr}x^v\partial_r.
$$
For any given $u $ satisfying $i'\in
u,$ since $|u|\geq 3,$ there is $l\in u$ such that $l\neq i'.$
Direct computation shows that
$[\mathrm{T_H}(x^{(a\varepsilon_i)}x_{i'}),\,\Delta_{l'}]=0.$
Applying $\phi$ yields
$$-\sum_{i'\in u\in
\mathbb{B}^1}c_{ui'}x^u\partial_{i'}+ \sum_{i'\not\in v\in
\mathbb{B}}(c_{vl}x^v\partial_l-c_{vl'}v^v\partial_{l'})-\delta_{l\in
v}\sum_{i'\not\in v\in \mathbb{B}\atop r\in Y\setminus
\{i,i'\}}c_{vr}x^v\partial_r=0.$$
 It follows that $c_{ui'}=0$ and therefore,
 $$
 \phi(\mathrm{T_H}(x^{(a\varepsilon_i)}x_{i'}))=\sum_{i'\not\in
v\in \mathbb{B}\atop r\in Y\setminus \{i,i'\}}c_{vr}x^v\partial_r.
$$
For any given $v$ satisfying $i'\not\in v,$ since $|v|\geq 3$, there is
$q\in v$
such that
$$[\mathrm{T_H}(x^{(a\varepsilon_i)}x_{i'}),\,\Delta_{q'}]=0.$$
Applying $\phi$ yields
\[\sum_{i'\not\in v\in \mathbb{B}}(c_{vq}x^v\partial_q-c_{vq'}x^v\partial_{q'})-\sum_{i'\not\in v\in
\mathbb{B}\atop r\in Y\setminus\{i,i'\}}c_{vr}x^v\partial_r=0.\] Therefore,
$c_{vr}=0$ for $ r\in Y\setminus q.$
Thus
\[
\phi(\mathrm{T_H}(x^{(a\varepsilon_i)}x_{i'}))=\sum_{i'\not\in
v\in \mathbb{B}^1}c_{vq}x^v\partial_q.
\]
Since $|v|\geq 3$,
there is $r\in v,\,q\neq r$ such that
$$[\mathrm{T_H}(x^{(a\varepsilon_i)}x_{i'}),\,\Delta_{r'}]=0.$$
Applying $\phi$, we have
$-\sum_{i'\not\in v\in \mathbb{B}^1}c_{vq}x^v\partial_q=0.$ It follows that $c_{vq}=0.$
So  \[\phi(\mathrm{T_H}(x^{(a\varepsilon_i)}x_{i'}))=0.\]

\noindent Step (ii):
$a\equiv 2\pmod{p}.$ By (\ref{he3.4}) and (\ref{he3.5}),  $c_{ur}=0$ for $r\in
Y\setminus i'.$
It follows that
\[\phi(\mathrm{T_H}(x^{(a\varepsilon_i)}x_{i'}))=\sum_{i'\in
u\in \mathbb{B}^1}c_{ui'}x^u\partial_{i'}+\sum_{i'\not\in v\in
\mathbb{B}^1}c_{vi'}x^v\partial_{i'}.\]
For any fixed $u $ satisfying $i'\in u,$
by  (\ref{he3.6}),  $c_{ui'}=0.$ Consequently,
 $$\phi(\mathrm{T_H}(x^{(a\varepsilon_i)}x_{i'}))=\sum_{i'\not\in
v\in \mathbb{B}^1}c_{vi'}x^v\partial_{i'}.$$
For any fixed $v$ satisfying
 $i'\not\in v,$ since $|v|\geq 3$, there is  $ l\in v$ such that
\[
[\mathrm{T_H}(x^{(a\varepsilon_i)}x_{i'}),\,\Delta_{l'}]=0.
\]
Applying $\phi$ we have
$-\sum_{i'\not\in v\in \mathbb{B}^1}c_{vi'}x^v\partial_{i'}=0.$
This implies that $c_{vi'}=0 $ and therefore,
 $
 \phi(\mathrm{T_H}(x^{(a\varepsilon_i)}x_{i'}))=0.
 $\\

 \noindent Step (iii):
$a\equiv -1\pmod{p}.$ By (\ref{he3.4}) and (\ref{he3.6}),
$c_{ur}=0$ for $ r\in
Y\setminus i.$
For any fixed $v $ satisfying $i'\not\in v,$
by (\ref{he3.5}), $c_{vi}=0 $ and consequently,
 $$
 \phi(\mathrm{T_H}(x^{(a\varepsilon_i)}x_{i'}))=\sum_{i'\in
u\in \mathbb{B}^0}c_{ui}x^u\partial_i.
$$
For any fixed $u$ satisfying $i'\in u,$
since $|u|\geq 3,$ there is $l\in u$ such that $l\neq i'.$
Note that
$[\mathrm{T_H}(x^{(a\varepsilon_i)}x_{i'}),\,\Delta_{l'}]=0.$
Applying $\phi$ we obtain that
$-\sum_{i'\in u\in \mathbb{B}^0}c_{ui}x^u\partial_i=0.$ This implies that $c_{ui}=0 $
and hence, $\phi(\mathrm{T_H}(x^{(a\varepsilon_i)}x_{i'}))=0.$\\

\noindent Step (iv):
$a\not\equiv -1,1, 2\pmod{p}.$ By  (\ref{he3.4})--(\ref{he3.6}),
we have $c_{ur}=0$ for $ r\in
Y$ and $u\in \mathbb{B}.$
It follows that $\phi(\mathrm{T_H}(x^{(a\varepsilon_i)}x_{i'}))=0.$
 Note that
\[
[\mathrm{T_H}(x^{(a\varepsilon_i)}x_{i'}),\,\mathrm{T_H}(x_ix_k)]
=-\mathrm{T_H}(x^{(a\varepsilon_i)}x_k) \quad\mbox{for}\,
k\in Y_{}\setminus i'.
\]
Applying $\phi$ the equation above gives
$$\phi(\mathrm{T_H}(x^{(a\varepsilon_i)}x_k))=0\quad\mbox{for}\,
k\in Y_1\setminus i'.$$
Hence, $\phi(\mathrm{T_H}(x^{(a\varepsilon_i)}x_k))=0\quad\mbox{for}\, i\in
Y_0,\,k\in Y_1.$\end{proof}

We next consider the subset $N$ of the generating set of
$\mathcal{HO}.$ First consider the action on $N$ of derivations of
odd $\mathbb{Z}$-degree.

\begin{lemma}\label{ht3.2} Let $\phi\in
\mathrm{Der}_t(\mathcal{HO},\mathcal{W})$ with $t\geq 0 $  odd.
If $\phi(\mathcal{HO}_{-1}\oplus\mathcal{HO}_0)=0,$ then
$\phi(\mathrm{T_H}(x_kx_lx_q))=0 $ for all $k,l,q\in Y_1.$
\end{lemma}
\begin{proof} By Lemma \ref{ht1.2},
one may assume that
\[\phi(\mathrm{T_H}(x_kx_lx_q))=\sum_{r\in Y_1, u\in
\mathbb{B}^1}c_{ur}x^u\partial_r\quad\mbox{where}\, c_{ur}\in \mathbb{F}.\]
Note that
$[\mathrm{T_H}(x_kx_lx_q),\,\Delta_{k'}]=-\mathrm{T_H}(x_kx_lx_q).$
Applying $\phi$ gives
\[\sum_{u\in
\mathbb{B}^1}c_{uk}x^u\partial_k-\delta_{k\in u}\sum_{r\in Y_1,\,u\in
\mathbb{B}^1}c_{ur}x^u\partial_r=-\sum_{r\in Y_1,\,u\in
\mathbb{B}^1}c_{ur}x^u\partial_r.
\]
For any fixed $u,$ if $k\not\in u,$
then $c_{ur}=0$ for all $ r\in Y_1.$  If  $k\in u,$ then $c_{uk}=0.$
Therefore,
\[\phi(\mathrm{T_H}(x_kx_lx_q))=\sum_{k\in u\in
\mathbb{B}^1\atop r\in Y_1\setminus k}c_{ur}x^u\partial_r .\]
Similarly,  it is easily seen that
\[\phi(\mathrm{T_H}(x_kx_lx_q))=\sum_{k,l,q\in u\in
\mathbb{B}^1\atop r\in Y_1\setminus \{k,l,q\}}c_{ur}x^u\partial_r.
\]
The equation above implies that our conclusion holds when $n=3.$
Let us consider the case that $n>3.$ For $n=4,$ given any $u $ satisfying
 $k,\,l,\,q\in u,$ since $u\in \mathbb{B}^1 $ implies that $|u|<4,$ there is
 $
s\in Y_1\setminus\{k,l,q\} $ such that $s\not\in u.$ Direct computation shows that
 $[\mathrm{T_H}(x_kx_lx_q),\,\Delta_{s'}]=0.$
Applying $\phi$ yields
$\sum_{k,l,q\in u\in \mathbb{B}^1}c_{us}x^u\partial_s=0.$
Consequently,
$c_{us}=0 $ and therefore, $\phi(\mathrm{T_H}(x_kx_lx_q))=0.$ For $n>4,$
just as above, one may obtain that
 $$\sum_{k,l,q\in u\in \mathbb{B}^1}c_{us}x^u\partial_s-
\delta_{s\in u}\sum_{k,l,q\in u\in \mathbb{B}^1\atop r\in
Y_1\setminus \{k,l,q\}}c_{ur}x^u\partial_r=0.$$
 If $s\not\in u $ then $c_{us}=0 $
for all $ s\in Y_1\setminus\{k,l,q\};$ if  $s\in u  $ then $c_{ur}=0 $ for all
$r\in Y_1\setminus s.$
Consequently,
\[\phi(\mathrm{T_H}(x_kx_lx_q))=\sum_{k,l,q,s\in u\in
\mathbb{B}^1}c_{us}x^u\partial_s.\]
Note that
$$[\mathrm{T_H}(x_kx_lx_q),\,\mathrm{T_H}(x_sx_{r'})]=0\quad \mbox{for}\;\, r\in
Y_1\setminus\{k,l,q,s\}.$$
Applying $\phi$ yields
$\sum_{k,l,q,s\in u\in \mathbb{B}^1}c_{us}x^u\partial_r=0.$ This implies that $c_{us}=0 $ and
hence
 $$\phi(\mathrm{T_H}(x_kx_lx_q))=0 \quad\mbox{for all}\,\; k,l,q\in Y_1.$$
 The proof is complete. \end{proof}

Let us consider the derivations of even $\mathbb{Z}$-degree.

\begin{lemma}\label{ht3.3}Let $\phi\in
\mathrm{Der}_t(\mathcal{HO},\mathcal{W})$ with $t\geq 0 $  even.
If $\phi(\mathcal{HO}_0\oplus\mathcal{HO}_{-1})=0,$
then $\phi(\mathrm{T_H}(x_kx_lx_q))=0$ for all $k,l,q\in Y_1.$
\end{lemma}
\begin{proof} By Lemma \ref{ht1.2},
we may assume that
\[\phi(\mathrm{T_H}(x_kx_lx_q))=\sum_{r\in Y_0,\, u\in
\mathbb{B}^0}c_{ur}x^u\partial_r\quad\mbox{where}\,c_{ur}\in \mathbb{F}.
\]
Assume that $n>3.$ Then for arbitrary $s\in Y_1\setminus \{k,l,q\},$ it is easy to see that
$[\mathrm{T_H}(x_kx_lx_q),\,\Delta_{s'}]=0.$
Applying $\phi$ we have
\[
-\sum_{u\in \mathbb{B}^0}c_{us'}x^u\partial_{s'}-\delta_{s\in u}\sum_{r\in
Y_0,\,u\in \mathbb{B}^0}c_{ur}x^u\partial_r=0.\]
Hence,
\[\phi(\mathrm{T_H}(x_kx_lx_q))=\sum_{s\not\in u\in
\mathbb{B}^0\atop r\in Y_0\setminus s'}c_{ur}x^u\partial_r.\]
Note that
$[\mathrm{T_H}(x_kx_lx_q),\,\mathrm{T_H}(x_{k'}x_l)]=0.$
Applying $\phi$ we have
 $$
-\sum_{s\not\in u\in
\mathbb{B}^0}c_{u{k'}}x^u\partial_{l'}-\sum_{s\not\in u\in
\mathbb{B}^0\atop r\in Y_0\setminus s'}c_{ur}x_l\partial_k(x^u)\partial_r=0.
$$
It follows that
$c_{uk'}=0.$
 Similarly, one  gets $c_{ul'}=0 $ and  $c_{uq'}=0.$
 It follows that $\phi(\mathrm{T_H}(x_kx_lx_q))=0.$
 For $n=3,$ the argument is similar and much easier.
\end{proof}

By Lemmas \ref{ht3.1}--\ref{ht3.3}, we have the following proposition.

\begin{proposition}\label{ht3.4} Let $\phi\in
\mathrm{Der}_t(\mathcal{HO},\mathcal{W})$ with $t\geq 0.$
If $\phi(\mathcal{HO}_{-1}\oplus\mathcal{HO}_0)=0,$ then $\phi=0$.
\end{proposition}

This proposition tells us that
the nonnegative $\mathbb{Z}$-degree derivations from
$\mathcal{HO}$ into $\mathcal{W}$ are completely determined by the top of $\mathcal{HO}$.
In order to reduce the nonnegative  $\mathbb{Z}$-degree derivations from
$\mathcal{HO}$ into $\mathcal{W}$ to vanish on the top, we establish the following lemma.

\begin{lemma}\label{ht3.5} Let $r\leq n$ be a positive integer and
$f_1,\,f_2,\,\ldots,\,f_r\in \Lambda(n).$
 Suppose that
\begin{enumerate}
\item [$\mathrm{(a)}$]\,\,
$ \Delta_{i} (f_j)=\Delta_{j}(f_i),\quad 1\leq
i,j\leq r;$

\item [$\mathrm{(b)}$]\,\,$\Delta_{i}(f_i)=f_i,\,i=1,\ldots,r.$
\end{enumerate}

Then there is  $f\in \Lambda(n) $
such that $\Delta_{i}(f)=f_i $ for $i=1,\ldots,r.$
\end{lemma}

\begin{proof}
We verify the conditions of Lemma \ref{ht1.3}. Since
$(x_{k}\partial_{k})^{2}=x_{k}\partial_{k} $ for all $k\in Y_{1},$
 it is easily seen that
$\Delta_{i}^2|_{\Lambda(n)}=\Delta_{i}|_{\Lambda(n)}.$
Let $A_i:=\Delta_{i}|_{\Lambda(n)},\,B_i:=\mathrm{id}_{\Lambda{(n)}}.$
Then $A_iB_iA_i=A_i^2=A_i,$ that is, (i)  holds.
Lemma \ref{ht1.3}(ii) is obvious. By (a), (iii) holds. By  (b),
$A_iB_i(f_i)=A_i(f_i)=f_i.$ Clearly,  $ A_iB_j=B_jA_i,\quad 1\leq i,j\leq
r,\,i\neq j.$ Thus (iv) holds. By Lemma \ref{ht1.3}, there is $f\in
\Lambda(n)$ such that $\Delta_{i}(f)=f_i $ for
$i=1,\ldots,r.$\end{proof}

According to Proposition \ref{ht3.4}, it suffices to
consider the top for the nonnegative $\mathbb{Z}$-degree derivations from
 $\mathcal{HO}$ into $\mathcal{W}$. In the proof of the following lemma we
 shall adopt the methods used in the proof of
 \cite[Lemma 4.2.5, Proposition 3.2.4]{lz1}.

\begin{lemma}\label{ht3.6} Let $\phi\in
\mathrm{Der}_t(\mathcal{HO},\mathcal{W})$ with $t\geq 0 $  even.
If $\phi(\mathcal{HO}_{-1})=0,$ then there is $D\in \mathcal{G}_t $
such that $(\phi-\mathrm{ad}D)(\mathcal{HO}_{-1}\oplus\mathcal{HO}_0)=0.$
\end{lemma}

\begin{proof} By Lemma \ref{ht1.2},
 one may assume that
 $\phi(\Delta_{i})=\sum_{r\in Y_1}f_{ri'}\partial_r,$
where $f_{ri'}\in \Lambda(n).$ Note that
$$
[\Delta_{i},\Delta_{j}]=0 \quad\mbox{for all}\,\;
i,j \in Y_0\,\mbox{with}\,\; i\not=j.
$$
Applying $\phi$  and comparing  coefficients one may obtain that
\begin{equation}\label{he3.7}
\Delta_i(f_{rj'})=\Delta_j(f_{ri'}) \quad\mbox{for}\; \,
r\in Y_1\setminus\{i',j'\};
\end{equation}
\begin{equation}\label{he3.8}
\Delta_j(f_{i'i'})=\Delta_i(f_{i'j'})-f_{i'j'} \quad\mbox{whenever}\;\, i,j \in Y_0\,\,\mbox{with}\,\; i\not=j.
\end{equation}
For $r,\,i'\in Y_1,$
Suppose $ f_{ri'}=\sum_{|u|=t+1}c_{uri'}x^u,$  where
$c_{uri'}\in \mathbb{F}.$ Then by (\ref{he3.7}),
$\delta_{i'\in u}c_{urj'}=\delta_{j'\in u}c_{uri'}$ for $ r\neq
i',j'.$ This implies that
$$c_{uri'}\neq 0\;\,\mbox{and}\,\;j'\in u\Longleftrightarrow
c_{urj'}\neq 0\;\,\mbox{and}\;\,i'\in u.$$
Let ${r,i'}\in Y_1$ with $r\neq i'.$
If $c_{uri'}\neq 0,$ then the implication ensures that there is $i'\in u $ and
therefore,
\begin{equation}\label{he3.9}
\Delta_i(f_{ri'})=f_{ri'},\quad r\neq
i'.
\end{equation}
 For any fixed  $r\in Y_1,$
by Lemma \ref{ht3.5}, there is $\overline{f}_r\in \Lambda(n) $ such that
\begin{equation}\label{he3.10}
\Delta_i(\overline{f}_r)=f_{ri'} \quad\mbox{for all }\,
i'\in Y_1\setminus r.
\end{equation}

Assert that $\Delta_i(f_{i'i'})=0 $ for all
$i'\in Y_1.$
We first consider the case $t\geq 2.$ By (\ref{he3.8}),
$$
 \Delta_j\Delta_i
(f_{i'i'}) = \Delta_i\Delta_j(f_{i'i'})
=\Delta_i^2(f_{i'j'})-\Delta_i(f_{i'j'}) =0 \quad\mbox{for }\,
 j'\in Y_1\setminus{i'}.
$$
Assume that $\Delta_i(f_{i'i'})\neq 0.$
Since $\mathrm{zd}(f_{i'i'})=t+1\geq 3,$ there is $r\in u\setminus {i'} $
such that $ \Delta_{r'} \Delta_i(f_{i'i'})\neq 0,$
contradicting the equation above. Hence, $\Delta_i(f_{i'i'})=0.$

 Let us consider the case  $t=0.$ Apply $\phi$ to the equation that
 $[\mathrm{T_H}(x_ix_{j'}),\Delta_{i}]=\mathrm{T_H}(x_ix_{j'})$
for $j\in Y_0\setminus i.$  We have
\[\phi(\mathrm{T_H}(x_ix_{j'}))+[\Delta_i,\phi(\mathrm{T_H}(x_ix_{j'}))]
=[\mathrm{T_H}(x_ix_{j'}),\phi(\Delta_i)].\]
On the other hand, since $\mathrm{zd}(f_{ri'})=1,$ by (\ref{he3.9}),
for $r\neq i',$ we have $f_{ri'}=c_{ri'}x_{i'},$ where $ c_{ri'}\in
\mathbb{F}.$
Then
\begin{equation*}\phi(\Delta_i)=\sum_{r\in
Y_1}f_{ri'}\partial_r=f_{i'i'}\partial_{i'}+\sum_{r\in Y_1\setminus
i'}c_{ri'}x_{i'}\partial_r.
\end{equation*}
Consequently,
\begin{equation}
\phi(\mathrm{T_H}(x_ix_{j'})) + [\Delta_i,\,\phi(\mathrm{T_H}(x_ix_{j'}))]
=(\partial_{i'}(f_{i'i'})x_{j'}-c_{j'i'}x_{i'})\partial_{i'}+\sum_{r\in
Y_1\setminus{i'}}c_{ri'}x_{j'}\partial_r.\label{he3.11}
\end{equation}
Since $\mathrm{zd}(\phi)=0,$
 $\phi(\mathrm{T_H}(x_ix_{j'}))=\sum_{s,r\in
Y_1}\mu_{s,r}^{(i,j')}x_s\partial_r \,\,\mbox{where}\,\mu_{s,r}^{(i,j')}\in
\mathbb{F}.$ Note that
$
[\Delta_i,x_s\partial_r]=(\delta_{i's}-\delta_{i'r})x_s\partial_r.$
If follows that
\[\phi(\mathrm{T_H}(x_ix_{j'}))+[\Delta_i,\,
\phi(\mathrm{T_H}(x_ix_{j'}))]=\sum_{s,r\in
Y_1}\mu_{s,r}^{(i,j')}x_s\partial_r +\sum_{s,r\in
Y_1}(\delta_{i's}-\delta_{i'r})\mu_{s,r}^{(i,j')}x_s\partial_r.\]
In the equation above the coefficient of $\partial_{i'}$ is
\[\sum_{s\in
Y_1}\mu_{s,i'}^{(i,j')}x_s+\sum_{s\in
Y_1}(\delta_{i's}-\delta_{i'i'})\mu_{s,i'}^{(i,j')}x_s=\sum_{s\in
Y_1}\delta_{i's}\mu_{s,i'}^{(i,j')}x_s=\mu_{i',i'}^{(i,j')}x_{i'}.\]
By  (\ref{he3.11}),
$\partial_{i'}(f_{i'i'})x_{j'}-c_{j'i'}x_{i'}=\mu_{i',i'}^{(i,j')}x_{i'}$ for $
i\neq j.$ Hence, $\partial_{i'}(f_{i'i'})=0 $ and
 $\Delta_i(f_{i'i'})=0.$ Thus our assertion holds.

 For $r\in Y_1,$
let $f_r:=-f_{rr}+ \Delta_{r'} (\overline{f}_r).$
Clearly, $f_r\in \Lambda(n) $
and
$$\Delta_{r'}(f_r)=\Delta_{r'}(-f_{rr})+
\Delta_{r'}^2(\overline{f}_r)
=\Delta_{r'}(\overline{f}_r)=f_r+f_{rr}.$$
Then
\begin{equation}\label{he3.12}
\Delta_{r'}(f_r)-f_r=f_{rr}.
\end{equation}
For  $i'\in Y_1\setminus r,$
by (\ref{he3.8}) and (\ref{he3.10}), one may compute
$\Delta_i(f_r)
 = f_{ri'}.
$
Putting $D':=-\sum_{r\in Y_1}f_r\partial_r,$ then by   (\ref{he3.12}), we have
$$
[D',\Delta_i] = -\sum_{r\in
Y_1}[f_r\partial_r,\Delta_i]  =\sum_{r\in
Y_1\setminus
i'}f_{ri'}\partial_r+f_{i'i'}\partial_{i'}=\phi(\Delta_i).
$$
Put $D :=D_t'.$ Then $D\in \mathcal{G}_t $ and
$[D,\,\Delta_i]=\phi(\Delta_i).$ It follows that
$(\phi-\mathrm{ad}D)(\Delta_{i})=0.$
Let $\psi:=\phi-\mathrm{ad}D $ and suppose
 \[\psi(\mathrm{T_H}(x_ix_k))=\sum_{r\in Y_1,\, u\in
\mathbb{B}^1}c_{ur}x^u\partial_r \quad\mbox{for}\;\, k\in Y_{1}\setminus i'.\]
Note that
$[\mathrm{T_H}(x_ix_k),\Delta_{i}]=\mathrm{T_H}(x_ix_k).$
Applying $\psi$ we have
\[\sum_{u\in
\mathbb{B}^1}c_{ui'}x^u\partial_{i'}-\delta_{i'\in u}\sum_{r\in Y_1,\,u\in
\mathbb{B}^1}c_{ur}x^u\partial_r=\sum_{r\in Y_1,\,u\in
\mathbb{B}^1}c_{ur}x^u\partial_r.\]
Therefore,
$\psi(\mathrm{T_H}(x_ix_k))=\sum_{i'\not\in u\in
\mathbb{B}^1}c_{ui'}x^u\partial_{i'}.$ Note that
$$[\mathrm{T_H}(x_ix_k),\,\Delta_{k'})]
=-\mathrm{T_H}(x_ix_k)\quad\mbox{for}\;\ k\in Y_1\setminus i'.$$
Applying $\psi$ we have
 $$-\delta_{k\in
u}\sum_{i'\not\in u\in
\mathbb{B}^1}c_{ui'}x^u\partial_{i'}=-\sum_{i'\not\in u\in
\mathbb{B}^1}c_{ui'}x^u\partial_{i'}.
$$
Hence
$$\psi(\mathrm{T_H}(x_ix_k))=\sum_{k\in u\atop i'\not\in u\in
\mathbb{B}^1}c_{ui'}x^u\partial_{i'}.
$$
If $t>0,$ for any fixed $u $ satisfying
 $i'\not\in u$ and $k\in u,$ since $|u|\geq 2,$ there is $r\in u\setminus
k$ such that  $[\mathrm{T_H}(x_ix_k),\,\Delta_{r'}]=0.$
Then by applying  $\psi$ we have
$-\sum_{k\in u\atop i'\not\in u\in
\mathbb{B}^1}c_{ui'}x^u\partial_{i'}=0 $ and hence $c_{ui'}=0.$
Thus   $\psi(\mathrm{T_H}(x_ix_k))=0.$  If $t=0,$
then $\psi(\mathrm{T_H}(x_ix_k))=c_{ki'}x_k\partial_{i'}.$
Applying $\psi$ to the equation that
$$[\mathrm{T_H}(x_ix_k),\,\mathrm{T_H}(x_{i'}x_{l'})]=\mathrm{T_H}(x_{l'}x_k)\quad
\mbox{for}\;\, l\in Y_1\setminus \{k,i'\},$$ one may get
$c_{ki'}+c_{i'l}=c_{kl}.$ Similarly,
from the equation that
$$[\mathrm{T_H}(x_ix_k),\,\mathrm{T_H}(x_{k'}x_{i'})]
=\Delta_{k'}-\Delta_{i },$$
one   gets
  $c_{ki'}+c_{i'k}=0.$
Clearly, the following system of $n-1$ equations has solutions:
$$\lambda_{1'}-\lambda_{2'}=c_{1'2'}$$
$$\lambda_{1'}-\lambda_{3'}=c_{1'3'}$$
$$\cdots\cdots\cdots\cdots\cdots$$
$$\lambda_{1'}-\lambda_{n'}=c_{1'n'}.$$
Let $(\lambda_{1'},\lambda_{2'},\dots,\lambda_{n'})$ be a solution. Then
\begin{equation*}\lambda_k-\lambda_l=(\lambda_k-\lambda_{1'})+(\lambda_{1'}-\lambda_l)
=c_{kl}.
\end{equation*}
Let $D'':=\sum_{r\in Y_1}\lambda_rx_r\partial_r$ and
 $\varphi:=\psi-\mathrm{ad}D''.$
 Then
$\varphi(\mathrm{T_H}(x_ix_k))=0
$ for $ i\not=k'.$
In addition,
$
\varphi(\Delta_{i})
=0.
$
The proof is complete.\end{proof}

\begin{lemma}\label{ht3.7} Let $\phi\in
\mathrm{Der}_t(\mathcal{HO},\mathcal{W})$ with $t\geq 0$  odd.
If $\phi(\mathcal{HO}_{-1})=0,$ then there is $D\in \mathcal{G}_t $ such that
 $(\phi-\mathrm{ad}D)(\mathcal{HO}_{-1}\oplus\mathcal{HO}_0)=0.$
\end{lemma}

\begin{proof} By Lemma \ref{ht1.2}, assume that
 $\phi(\Delta_{i})=\sum_{r\in Y_0}f_{ri}\partial_r,$
where $f_{ri}\in \Lambda(n).$ For arbitrary $j\in Y_0\setminus i,$ we have
 $[\Delta_{i},\,\Delta_{j}]=0.$
Then
\[-f_{ji}\partial_j-\sum_{r\in Y_0}\Delta_j(f_{ri})\partial_r+\sum_{r\in
Y_0}\Delta_i(f_{rj})\partial_r+f_{ij}\partial_i=0.\]
A comparison of coefficients shows that
\begin{equation}
\Delta_j(f_{ri})=\Delta_i(f_{rj})\quad\mbox{for}\,\;
r\in Y_{0}\setminus \{ i,j\},\label{he3.13}
\end{equation}
\begin{equation}
 \Delta_j(f_{ii})=\Delta_i(f_{ij})+f_{ij} \quad\mbox{for}\,\;i,j\in Y_0\,\;\mbox{with}\,\;
i\neq j.\label{he3.14}
\end{equation}
Suppose   $f_{ri}=\sum_{|u|=t+1}c_{uri}x^u,$ where $
c_{uri}\in \mathbb{F}.$ By (\ref{he3.13}),
$\delta_{j'\in
u}c_{uri}=\delta_{i'\in u}c_{urj}$ for $ r\in Y_{0}\setminus \{ i,j\}.$ Then we have the following
implication:
$$j'\in
u\,\mbox{and }\,c_{uri}\neq 0\Longleftrightarrow i'\in
u\,\;\mbox{and}\,\; c_{urj}\neq 0, \;\, r\in Y_{0}\setminus \{ i,j\}.$$
Let $r,i\in Y_0$ with $r\neq i.$
If $c_{uri}\neq 0,$ then $i'\in u $ and therefore,
 $\Delta_i(f_{ri})=f_{ri}.$
 For any fixed $r\in
Y_0$, by Lemma  \ref{ht3.5}, there is $ {f}_r\in
\Lambda(n) $ such that
\begin{equation}
\Delta_i( {f}_r)=f_{ri} \quad\mbox{for all}\,\;
i\in Y_0\setminus r.\label{he3.15}
\end{equation}
 For any fixed $r\in Y_0$
and  $i\in Y_{0}\setminus   r$, by (\ref{he3.14}) and (\ref{he3.15}), we have
\begin{eqnarray*}
  \Delta_i(f_{rr}-f_r- \Delta_r  (f_r))
  &=&\Delta_i(f_{rr})-\Delta_i(f_r)-\Delta_i \Delta_r  (f_r)\\
  &=&\Delta_i(f_{rr})-f_{ri}-\Delta_i \Delta_r  (f_r)\\
  &=&\Delta_r  (f_{ri})- \Delta_r  \Delta_i(f_r)\\
  &=&\Delta_r  \Delta_i({f}_r)
- \Delta_r  \Delta_i(f_r)\\
&=&0.
\end{eqnarray*}
Let
$
f:=f_{rr}-f_r- \Delta_r  (f_r)=\sum_vc_vx^v.
$
Assume that $f\neq 0.$ Then there is $v $ such that $c_v\neq 0.$ Since $|v|\geq 2,$
there is $i'\in v\setminus r'.$
Then $0=\Delta_i(f)=c_vx^v $ and therefore, $c_v=0,  $ a contradiction.
This shows that $f=0,$ that is,
\begin{equation}
f_{rr}= \Delta_r  (f_r)+f_r \quad \mbox{for}\;\, r\in
Y_0.\label{he3.17}
\end{equation}
Let $D':=-\sum_{r\in Y_0}f_r\partial_r.$ By
  (\ref{he3.15}) and (\ref{he3.17}),
  \begin{eqnarray*}
[D',\Delta_i]&=&-\sum_{r\in
Y_0}[f_r\partial_r,\,\Delta_i]\\&=&\sum_{r\in
Y_0}\Delta_i(f_r)\partial_r+f_i\partial_i\\&=&\sum_{r\in
Y_0\setminus
i}\Delta_i(f_r)\partial_r+(\Delta_i(f_i)+f_i)\partial_i\\
&=&\sum_{r\in Y_0\setminus i}f_{ri}\partial_r+f_{ii}\partial_i\\
&=&\phi(\Delta_i).\end{eqnarray*}
Let $D:=D'_t.$ Then $D\in
\mathcal{G}_t $ and
$[D,\,\Delta_i]=\phi(\Delta_i).$
Put $\psi:=\phi-\mathrm{ad}D.$  Then $\psi(\Delta_i)=0.$
Assume that \[\psi(\mathrm{T_H}(x_ix_k))=\sum_{r\in Y_0,\,u\in
\mathbb{B}^0}c_{ur}x^u\partial_r \quad\mbox{where}\,\; c_{ur}\in\mathbb{F},\,k\in Y_{1}\setminus
 i'.\]
By applying $\psi$ to
$[\mathrm{T_H}(x_ix_k),\Delta_i]=\mathrm{T_H}(x_ix_k),$
one  computes
\[-\sum_{u\in \mathbb{B}^0}c_{ui}x^u\partial_i-\delta_{i'\in u}\sum_{r\in Y_0,\,u\in
\mathbb{B}^0}c_{ur}x^u\partial_r=\sum_{r\in Y_0,\,u\in
\mathbb{B}^0}c_{ur}x^u\partial_r.\]
It follows that $c_{ur}=0 $ for all $r\in Y_0.$ Therefore,
$\psi(\mathrm{T_H}(x_kx_i))=0. $ Now we have proved that
 $(\phi-\mathrm{ad}D)(\mathcal{HO}_{-1}\oplus\mathcal{HO}_0)=0.$
\end{proof}

The following proposition asserts that all the derivations of nonnegative $\mathbb{Z}$-degree are inner.
\begin{proposition}\label{ht3.8}
$\mathrm{Der}_t(\mathcal{HO},\mathcal{W})=\mathrm{ad}\mathcal{W}_t$ for $
t\geq 0.$
\end{proposition}

\begin{proof} Clearly, $\mathrm{ad}\mathcal{W}_t\subset
\mathrm{Der}_t(\mathcal{HO},\mathcal{W}).$
Let $\phi\in \mathrm{Der}_t(\mathcal{HO},\mathcal{W}).$
By Lemma \ref{ht1.4}, there is $E\in \mathcal{W}_t $ such that
 $(\phi-\mathrm{ad}E)(\mathcal{HO}_{-1})=0.$ By Lemmas \ref{ht3.6} and
\ref{ht3.7}, there is $D\in \mathcal{G}_t $ such that
 $(\phi-\mathrm{ad}E-\mathrm{ad}D)(\mathcal{HO}_{-1}\oplus\mathcal{HO}_0)=0.$
By Proposition \ref{ht3.4}, $\phi-\mathrm{ad}E-\mathrm{ad}D=0 .$  Hence
 $\phi=\mathrm{ad}E+\mathrm{ad}D\in
 \mathrm{ad}\mathcal{W}_t.$\end{proof}

\section{Negative $\mathbb{Z}$-degree derivations from  $\mathcal{HO}$ into
$\mathcal{W}$}

In this section, we first prove that $\mathbb{Z}$-degree $-1$ derivations from
$\mathcal{HO}$ into $\mathcal{W}$ vanishing on  $\mathcal{HO}_0$ are necessarily zero,
then determine the  $\mathbb{Z}$-degree $-1$ derivations. For our purpose, we need the following two lemmas.

\begin{lemma}\label{ht4.1} Let $\phi\in
\mathrm{Der}_{-1}(\mathcal{HO},\mathcal{W}) $
satisfy $\phi(\mathcal{HO}_0)=0.$ Then $\phi(\mathrm{T_H}(x_kx_lx_q))=0 $
for all $ k,\,l,\,q\in Y_1.$
\end{lemma}
\begin{proof} By Lemma \ref{ht1.2},
assume that $\phi(\mathrm{T_H}(x_kx_lx_q))=\sum_{s,r\in Y_1}c_{sr}x_s\partial_r,$
where $c_{sr}\in \mathbb{F}.$  Note that
$[\Delta_{k'},\mathrm{T_H}(x_kx_lx_q)]=\mathrm{T_H}(x_kx_lx_q).$
Applying $\phi$ one gets
\[\sum_{r\in Y_1}c_{kr}x_k\partial_r-\sum_{s\in Y_1}c_{sk}x_s\partial_k=\sum_{s,r\in
Y_1}c_{sr}x_s\partial_r.\] A comparison of coefficients shows that
\[c_{kk}x_k-\sum_{s\in Y_1}c_{sk}x_s=\sum_{s\in Y_1}c_{sk}x_s,\]
and
\[c_{kr}x_k=\sum_{s\in Y_1}c_{sr}x_s,\quad r\in Y_1\setminus k.\]
Then we have $c_{sk}=0$ for $ s\in Y_1$ and $ c_{sr}=0$ for  $s,r\in Y_1\setminus
k.$ Thus, \[\phi(\mathrm{T_H}(x_kx_lx_q))=\sum_{r\in Y_1\setminus
k}c_{kr}x_k\partial_r.\] Applying $\phi$ to the equation that
$[\mathrm{T_H}(x_kx_lx_q),\mathrm{T_H}(x_{k'}x_l)]=0 $
 one may computes
\[c_{kl}x_k\partial_k-\sum_{r\in Y_1\setminus k}c_{kr}x_l\partial_r=0.\] It follows that
¹Ê$c_{kr}=0 $ for all $r\in Y_1\setminus
k.$ Hence $\phi(\mathrm{T_H}(x_kx_lx_q))=0.$
\end{proof}
\begin{lemma} \label{ht4.2} Suppose  $\phi\in
\mathrm{Der}_{-1}(\mathcal{HO},\mathcal{W}) $
 and $\phi(\mathcal{HO}_0)=0.$ Then $\phi(\mathrm{T_H}(x^{(a\varepsilon_
i)}x_k))=0 $ for all $i\in Y_0,\,k\in Y_1.$
\end{lemma}

\begin{proof} Proceed by induction on $a$. For $a=2$, by Lemma \ref{ht1.2}, one
may assume that
 $$\phi(\mathrm{T_H}(x^{(2\varepsilon_ i)}x_k))=\sum_{q,r\in
Y_1}c_{qr}x_q\partial_r \quad\mbox{where}\, c_{qr}\in \mathbb{F}.$$
Note that
$[\Delta_i,\,\mathrm{T_H}(x^{(2\varepsilon_
i)}x_k)]=(\delta_{i'k}-2)\mathrm{T_H}(x^{(2\varepsilon_i)}x_k).$
Applying $\phi$ and computing we have
\[\sum_{r\in Y_1}c_{i'r}x_{i'}\partial_r-\sum_{q\in Y_1}c_{qi'}x_q\partial_{i'}=(\delta_{i'k}-2)\sum_{q,r\in
Y_1}c_{qr}x_q\partial_r.\] It follows that
\[c_{i'i'}x_{i'}-\sum_{q\in Y_1}c_{qi'}x_q=(\delta_{i'k}-2)\sum_{q\in
Y_1}c_{qi'}x_q,\] and
\[c_{i'r}x_{i'}=(\delta_{i'k}-2)\sum_{q\in Y_1}c_{qr}x_q \quad\mbox{for}\;\, r\in Y_1\setminus i'.\]
Consequently, $c_{i'i'}=0$ and
$(\delta_{i'k}-1)c_{qi'}=0$ for $ q\in Y_1\setminus i';$ $ c_{qr}=0 $ for
$q\in Y_1,$ $ r\in
Y_1\setminus i'.$
Then
\[\phi(\mathrm{T_H}(x^{(2\varepsilon_i)}x_{k}))=\delta_{i'k}\sum_{q\in
Y_1\setminus i'}c_{qi'}x_q\partial_{i'}.
\]
For $l\in Y_1\setminus i',$
$[\mathrm{T_H}(x^{(2\varepsilon_i)}x_{i'}),\Delta_{l'}]=0.$
Applying $\phi$ and one may get
 $c_{li'}=0 $ and therefore,
  $\phi(\mathrm{T_H}(x^{(2\varepsilon_i)}x_{k}))=0.$  For $a\geq 3,$
by Lemma \ref{ht1.2},
 one may assume that
 \[\phi(\mathrm{T_H}(x^{(a\varepsilon_i)}x_k))=\sum_{r\in Y,u\in
\mathbb{B}}c_{ur}x^u\partial_r \quad \mbox{for all}\,\;i\in Y_0,\,k\in Y_1.\]
If $a$ is even, then
\[\phi(\mathrm{T_H}(x^{(a\varepsilon_i)}x_k))=\sum_{r\in
Y_1,\, u\in \mathbb{B}^1}c_{ur}x^u\partial_r \quad\mbox{for all}\,\;i\in
Y_0,\,k\in Y_1.\]
 For any fixed $u,$ since $|u|\geq 3,$ there is $j'\in u\setminus
\{i',k\}.$ Note that
$[\Delta_{j},\mathrm{T_H}(x^{(a\varepsilon_i)}x_k)]=0.$
Applying $\phi$ one gets
\[ \sum_{r\in Y_1,\,u\in \mathbb{B}^1}c_{ur}x^u\partial_r-\sum_{u\in
\mathbb{B}^1}c_{uj'}x^u\partial_{j'}=0.\] This implies that $c_{ur}=0$ for $r\in
Y_1\setminus j' $ and therefore,
 \[\phi(\mathrm{T_H}(x^{(a\varepsilon_i)}x_k))=\sum_{u\in
\mathbb{B}^1}c_{uj'}x^u\partial_{j'}.\]
Since $n\geq 3,$ take $l\in
Y_1\setminus \{k,j'\}.$
Applying $\phi$ to
$[\mathrm{T_H}(x^{(a\varepsilon_i)}x_k),\mathrm{T_H}(x_{j'}x_{l'})]=0,$ we have
$\sum_{u\in \mathbb{B}^1}c_{uj'}x^u\partial_l=0.$ This implies that $c_{uj'}=0 $ and hence
 $\phi(\mathrm{T_H}(x^{(a\varepsilon_i)}x_k))=0.$

If $a$ is odd, then  $\phi(\mathrm{T_H}(x^{(a\varepsilon_i)}x_k))=\sum_{r\in
Y_0,u\in \mathbb{B}^0, |u|\geq 2}c_{ur}x^u\partial_r.  $
Given  $u,$ if $|u|=2$ and $\{k,i'\}\not\subset u; $ or $|u|=2$ and $k=i';$ or
$|u|>2,$
then there is $j'\in u $ such that  $j'\neq k,\,i'.$
 Note that
$[\Delta_{j},\,\mathrm{T_H}(x^{(a\varepsilon_
i)}x_k)]=0.$
Applying $\phi$ we have
\[\sum_{r\in
Y_0,u\in \mathbb{B}^0}c_{ur}x^u\partial_r+\sum_{u\in
\mathbb{B}^0}c_{uj}x^u\partial_j=0.\] This yields  $c_{ur}=0$ for $r\in Y_0 .$
Thus
$$\phi(\mathrm{T_H}(x^{(a\varepsilon_ i)}x_k))=\sum_{r\in
Y_0}c_{vr}x^v\partial_r\quad\mbox{where}\;\,  v=\{i', k\}.$$
Clearly, it suffices to consider the case $a=3.$
Compute
$$[\Delta_{i},\,\mathrm{T_H}(x^{(3\varepsilon_
i)}x_k)]=-3\mathrm{T_H}(x^{(3\varepsilon_ i)}x_k).$$
Applying $\phi$ to the equation above we have
\[\sum_{r\in
Y_0}c_{vr}x^v\partial_r+c_{vi}x^v\partial_i=-3\sum_{r\in Y_0}c_{vr}x^v\partial_r,\]
 and consequently,  $ c_{vr}=0$
for $ r\in Y_0\setminus
i$ and $5c_{vi}=0.$
If $5\not\equiv 0\pmod{p},$ then $c_{vi}=0 $ and
 $$\phi(\mathrm{T_H}(x^{(3\varepsilon_ i)}x_k))=0.$$
 If $5\equiv
0\pmod{p},$ then
\[\phi(\mathrm{T_H}(x^{(3\varepsilon_ i)}x_k))=
 c_{vi}x^v\partial_i.\]
Applying to $[\mathrm{T_H}(x^{(3\varepsilon_ i)}x_k),
\mathrm{T_H}(x_{i}x_{k})]=0$  yields $c_{vi}=0 $ and therefore,
 $\phi(\mathrm{T_H}(x^{(3\varepsilon_i)}x_k))=0.$
\end{proof}

As a direct consequence of Lemmas \ref{ht4.1} and \ref{ht4.2}, we have the
following proposition.

\begin{proposition}\label{ht4.4}
$\mathrm{Der}_{-1}(\mathcal{HO},\mathcal{W})=\mathrm{ad}\mathcal{HO}_{-1}.$
\end{proposition}
\begin{proof} Let $\phi\in
\mathrm{Der}_{-1}(\mathcal{HO},\mathcal{W}).$ By Lemma \ref{ht1.2},
assume that $\phi(\mathrm{T_H}(x_ix_k))=\sum_{r\in Y_0}c_{ikr}\partial_r,$
where $c_{ikr}\in \mathbb{F},\,i\in Y_0,k\in Y_1.$
Applying $\phi$ to
 $[\mathrm{T_H}(x_ix_k),\Delta_{i}]=\mathrm{T_H}(x_ix_k)$ for $k\in
Y_1\setminus i',$  one may obtain that
$-c_{iki}\partial_i+c_{ii'i}\partial_{k'}=\sum_{r\in Y_0}c_{ikr}\partial_r,$ and therefore,
$$
c_{ikr}=0,\quad r\in Y_0\setminus k';\quad
  c_{ii'i}=c_{ikk'}.
$$ Thus,
$\phi(\mathrm{T_H}(x_ix_k))=c_{ikk'}\partial_{k'}=c_{ii'i}\partial_{k'}$ for
$
k\in Y_1\setminus i'.$ Let \[\psi:=\phi+\sum_{r\in
Y_0}c_{rr'r}\mathrm{ad}\partial_r\quad\mbox{where}\,\; c_{r
r'r}\in \mathbb{F}.\]
Then for $
k\in Y_1\setminus
i',$
\begin{eqnarray*}
\psi(\mathrm{T_H}(x_ix_k))&=&\phi(\mathrm{T_H}(x_ix_k))+\sum_{r\in
Y_0}c_{rr'r}[\partial_r,\,\mathrm{T_H}(x_ix_k)]\\&=&c_{ii'i}\partial_{k'}-c_{ii'i}\partial_{k'}\\&=&0.
\end{eqnarray*}
Note that
$[\Delta_{i},\,\Delta_{j}]=0$ for
$j\in Y_0\setminus i.$ Applying $\phi,$ one may obtain that
    $-c_{ii'j}\partial_j+c_{jj'i}\partial_i=0.$
  Hence, $c_{ii'j}=0 $ and  $\phi(\Delta_{i})=c_{ii'i}\partial_i.$ Thus
\begin{eqnarray*}\psi(\Delta_{i})
&=&\phi(\Delta_{i})+\sum_{r\in
Y_0}c_{rr'r}[\partial_r,\,\Delta_{i}]\\
&=&c_{ii'i}\partial_i-c_{ii'i}\partial_i\\&=&0.
\end{eqnarray*}
So far we have proved that $\psi(\mathcal{HO}_0)=0.$ By Lemmas \ref{ht4.1} and \ref{ht4.2},
 $\psi=0.$ Therefore,
 $$\phi=-\sum_{r\in Y_0}c_{rr'r}\mathrm{ad}\partial_r\in
\mathrm{ad}\mathcal{HO}_{-1},$$
proving that
$\mathrm{Der}_{-1}(\mathcal{HO},\mathcal{W})=\mathrm{ad}\mathcal{HO}_{-1}.$
\end{proof}

In order to determine the negative $\mathbb{Z}$-degree derivations
from $\mathcal{HO}$ into $\mathcal{W}$, we give the following lemma.

\begin{lemma}\label{ht4.5} Let $\phi\in
\mathrm{Der}_{-t}(\mathcal{HO},\mathcal{W})$ with $t>1.$ For $i\in
Y_0$ and $k\in Y_1,$ if $\phi(\mathrm{T_H}(x^{(t\varepsilon_ i)}x_k))=0,$
then $\phi=0.$
\end{lemma}
\begin{proof}
We first show that
$\phi(\mathrm{T_H}(x_kx_lx_q))=0$ for all $ k,l,q\in
Y_1.$ It suffices to consider the case  $\mathrm{zd}(\phi)=-2.$ By Lemma
\ref{ht1.2},
assume that
$$\phi(\mathrm{T_H}(x_kx_lx_q))=\sum_{r\in
Y_0}a_r\partial_r\quad \mbox{where}\,\; a_r\in \mathbb{F}.
$$
Note that
$$[\mathrm{T_H}(x_kx_lx_q),\,\Delta_{k'}]=-\mathrm{T_H}(x_kx_lx_q).$$
Applying $\phi$, one gets
$-a_{k'}\partial_{k'}=-\sum_{r\in Y_0}a_r\partial_r.$ This shows that $a_r=0$ for
 $
 r\in
Y_0\setminus k'.$ It follows that
 $$\phi(\mathrm{T_H}(x_kx_lx_q))=a_{k'}\partial_{k'}.$$ Similarly, noting that
$[\mathrm{T_H}(x_kx_lx_q),\,\mathrm{T_H}(x_{k'}x_l)]=0,$ one may show that
$a_{k'}=0 $ and therefore,
 $\phi(\mathrm{T_H}(x_kx_lx_q))=0.$

In the following we proceed by induction on $a$ to show that
  $\phi(\mathrm{T_H}(x^{(a\varepsilon_i)}x_k))=0 $  for all $i\in Y_0,k\in
Y_1.$ If $a\leq t$,
then $\phi(\mathrm{T_H}(x^{(a\varepsilon_i)}x_k))=0.$ Assume that $a>t.$
By inductive hypothesis and Lemma \ref{ht1.2}, one may assume that
\[\phi(\mathrm{T_H}(x^{(a\varepsilon_
i)}x_k))=\sum_{r\in Y,\
|u|=a-t}c_{ur}x^u\partial_r\quad\mbox{where}\,\; c_{ur}\in \mathbb{F}.\]

  Suppose
$a-t\geq 2.$ If $a-t$ is odd, then $a-t\geq 3 $ and
 \[\phi(\mathrm{T_H}(x^{(a\varepsilon_i)}x_k))=\sum_{r\in
Y_1,\, u\in \mathbb{B}^1}c_{ur}x^u\partial_r.\]
Just as in the proof of Lemma \ref{ht4.2} (for $a$ being even there),
one may show that $\phi(\mathrm{T_H}(x^{(a\varepsilon_
i)}x_k))=0.$ If $a-t$ is even,
then \[\phi(\mathrm{T_H}(x^{(a\varepsilon_
i)}x_k))=\sum_{r\in Y_0,u\in \mathbb{B}^0}c_{ur}x^u\partial_r.\]
Just as in the proof of Lemma \ref{ht4.2} (for $a$ odd there),
 one may show that $\phi(\mathrm{T_H}(x^{(a\varepsilon_
i)}x_k))=0.$

   Now consider the case $a-t<2.$ Note that in this case, we have  $a-t=1 $
 and  therefore,
  \[\phi(\mathrm{T_H}(x^{(a\varepsilon_ i)}x_k))=\sum_{q,r\in
Y_1}c_{qr}x_q\partial_r.\]

(1) Let $k\neq i'.$ Note that
$[\mathrm{T_H}(x^{(a\varepsilon_
i)}x_k),\,\Delta_{i}]=a\mathrm{T_H}(x^{(a\varepsilon_i)}x_k).$
Applying $\phi$ one gets
\[\sum_{q\in Y_1}c_{qi'}x_q\partial_{i'}-\sum_{r\in Y_1}c_{i'r}x_{i'}\partial_r=a\sum_{q,r\in
Y_1}c_{qr}x_q\partial_r.\] By comparing coefficients we have
\[\sum_{q\in Y_1}c_{qi'}x_q-c_{i'i'}x_{i'}=a\sum_{q\in
Y_1}c_{qi'}x_q,\]
and \[a\sum_{q\in Y_1}c_{qr}x_q+c_{i'r}x_{i'}=0\quad\mbox{for}\;\,
r\in Y_1\setminus i'.\] Consequently,
$$ac_{i'i'}=0;$$
$$(a-1)c_{qi'}=0\quad\mbox{for}\;\, q\in Y_1\setminus i';$$
$$(a+1)c_{i'r}=0\quad\mbox{for}\;\, r\in Y_1\setminus i';$$
$$ac_{qr}=0\quad\mbox{for}\;\, q,r\in Y_1\setminus i'.$$
If  $a\equiv 0\pmod{p},$
then $c_{qi'}=c_{i'r}=0$ for $q,r\in Y_1\setminus i' $ and it follows that
\[\phi(\mathrm{T_H}(x^{(a\varepsilon_
i)}x_k))=c_{i'i'}x_{i'}\partial_{i'}+\sum_{q,r\in Y_1\setminus
i'}c_{qr}x_q\partial_r.\] Note that
\[[\mathrm{T_H}(x^{(a\varepsilon_
i)}x_k),\,\Delta_{k'}]=-\mathrm{T_H}(x^{(a\varepsilon_i)}x_k).\]
Applying $\phi$ we have
\[\sum_{q\in Y_1\setminus
i'}c_{qk}x_q\partial_k-\sum_{r\in Y_1\setminus
i'}c_{kr}x_k\partial_r=-c_{i'i'}x_{i'}\partial_{i'}-\sum_{q,r\in Y_1\setminus
i'}c_{qr}x_q\partial_r.\]
Comparing coefficients yields
$$c_{i'i'}=0;$$
$$2\sum_{q\in Y_1\setminus i'}c_{qk}x_q-c_{kk}x_k=0;$$
$$\sum_{q\in Y_1\setminus i'}c_{qr}x_q-c_{kr}x_k=0\quad\mbox{for}\;\,
r\in Y_1\setminus \{i',k\}.$$ It follows that
$$ c_{i'i'}=0;\quad
  c_{qk}=0 \quad\mbox{for}\,\; q\in Y_1\setminus i'; \quad
 c_{qr}=0\quad\mbox{for}\,\;  q,r\in Y_1\setminus \{i',k\}.$$
Then
\[\phi(\mathrm{T_H}(x^{(a\varepsilon_i)}x_k))=\sum_{r\in
Y_1\setminus \{k,i'\}}c_{kr}x_k\partial_r.\] For $l\in Y_1\setminus
\{k,i'\},$ we have $[\mathrm{T_H}(x^{(a\varepsilon_
i)}x_k),\Delta_{l'}]=0.$
Applying $\phi$ yields  $c_{kl}=0.$ Thus,
$\phi(\mathrm{T_H}(x^{(a\varepsilon_ i)}x_k))=0.$

For the case with  $a\not\equiv 0\pmod{p},$ our discuss is divided into three parts.
First suppose  $a\equiv 1\pmod{p}.$
Then
$$c_{qr}=0 \quad \mbox{for}\, q,r\in Y_1\setminus i';\quad
 c_{i'r}=0\quad \mbox{for}\, r\in
Y_1.$$
Consequently,
\[\phi(\mathrm{T_H}(x^{(a\varepsilon_i)}x_k))=\sum_{q\in
Y_1\setminus i'}c_{qi'}x_q\partial_{i'}.\]
Note that
$[\mathrm{T_H}(x^{(a\varepsilon_
i)}x_k),\,\Delta_{k'}]=-\mathrm{T_H}(x^{(a\varepsilon_
i)}x_k).$
Applying $\phi$ one gets
\[-c_{ki'}x_k\partial_{i'}=-\sum_{q\in Y_1\setminus
i'}c_{qi'}x_q\partial_{i'}.\]
Hence $c_{qi'}=0$ for $ q\in Y_1\setminus
\{i',k\} $ and
 $$\phi(\mathrm{T_H}(x^{(a\varepsilon_i)}x_k))=c_{ki'}x_k\partial_{i'}.$$
For $l\in Y_1\setminus \{i',k\},$ it is easy to verify that
$[\mathrm{T_H}(x^{(a\varepsilon_
i)}x_k),\,\mathrm{T_H}(x_{i'}x_kx_l)]=0.$
Applying $\phi$ we have
$$[c_{ki'}x_k\partial_{i'},\,-x_kx_l\partial_i+x_{i'}x_l\partial_{k'}-x_{i'}x_k\partial_{l'}]+
[\mathrm{T_H}(x^{(a\varepsilon_i)}x_k),\,\phi(\mathrm{T_H}(x_{i'}x_kx_l))]=0.$$
Since $\phi(\mathrm{T_H}(x_{i'}x_kx_l))\in \mathcal{HO}_{-1},$
 one may assume that
 \[\phi(\mathrm{T_H}(x_{i'}x_kx_l))=\sum_{r\in
Y_0}a_r\partial_r \quad\mbox{where}\,\; a_r\in \mathbb{F}.\]
Then
\[c_{ki'}x_kx_l\partial_{k'}+\Big[\mathrm{T_H}(x^{(a\varepsilon_i)}x_k),\sum_{r\in
Y_0}a_r\partial_r\Big]=0.\]
Consequently,
\[c_{ki'}x_kx_l\partial_{k'}-a_ix^{((a-2)\varepsilon_i)}x_k\partial_{i'}
+a_ix^{((a-1)\varepsilon_i)}\partial_{k'}=0.\]
It follows that $c_{ki'}=0 $ and therefore,
 $\phi(\mathrm{T_H}(x^{(a\varepsilon_i)}x_k))=0.$
Second, suppose
$a\equiv -1\pmod{p}.$ Then $c_{qi'}=0$ for $ q\in Y_1 $ and $c_{qr}=0 $ for
$q,\,r\in Y_1\setminus i'.$
Thus, \[\phi(\mathrm{T_H}(x^{(a\varepsilon_i)}x_k))=\sum_{r\in
Y_1\setminus i'}c_{i'r}x_{i'}\partial_r.\] Note that
$[\mathrm{T_H}(x^{(a\varepsilon_
i)}x_k),\,\Delta_{k'}]=-\mathrm{T_H}(x^{(a\varepsilon_
i)}x_k).$
Applying $\phi$ one may verify that
  $c_{i'r}=0$ for $ r\in Y_1\setminus i'.$
Hence $$\phi(\mathrm{T_H}(x^{(a\varepsilon_i)}x_k))=0.$$
Third, suppose
$a\not\equiv -1,1\pmod{p}.$ Then it is clear that $\phi(\mathrm{T_H}(x^{(a\varepsilon_
i)}x_k))=0.$

 (2) Let us consider the case  $k=i'.$ Direct computation shows that
 $$
[\mathrm{T_H}(x^{(a\varepsilon_
i)}x_{i'}),\,\Delta_{i}]
=(a-1)\mathrm{T_H}(x^{(a\varepsilon_i)}x_{i'}).
$$
Applying $\phi$ and computing one may get
$$(a-2)\sum_{q\in
Y_1}c_{qi'}x_q+c_{i'i'}x_{i'}=0;$$ $$(a-1)\sum_{q\in
Y_1}c_{qr}x_q+c_{i'r}x_{i'}=0,\quad r\in Y_1\setminus i'.$$
Comparing coefficients yields that
$$(a-1)c_{i'i'}=0;$$ $$(a-2)c_{qi'}=0,\quad q\in Y_1\setminus
i';$$ $$ac_{i'r}=0,\quad r\in Y_1\setminus i';$$
$$(a-1)c_{qr}=0,\quad q,r\in Y_1\setminus
i'.$$
If $a\equiv 0\pmod{p},$
then
$\phi(\mathrm{T_H}(x^{(a\varepsilon_ i)}x_{i'}))=\sum_{r\in
Y_1\setminus i'}c_{i'r}x_{i'}\partial_r.$
For $ j\in Y_0\setminus i,$
it is easily verified that  $[\mathrm{T_H}(x^{(a\varepsilon_
i)}x_{i'}),\,\Delta_{j}]=0.$ Applying $\phi$ we have
\[\Big[\sum_{r\in Y_1\setminus
i'}c_{i'r}x_{i'}\partial_r,\,\Delta_{j}\Big]=0. \] Consequently,
$c_{i'j'}x_{i'}\partial_{j'}=0 $ and therefore,  ¹Ê$c_{i'j'}=0.$ It follows that
 $\phi(\mathrm{T_H}(x^{(a\varepsilon_ i)}x_{i'}))=0.$
If  $a\not\equiv 0\pmod{p},$ the following discussion is divided into three parts.\\

\noindent  Case (i):
$a\equiv 1\pmod{p}.$ Then
\[\phi(\mathrm{T_H}(x^{(a\varepsilon_
i)}x_{i'}))=c_{i'i'}x_{i'}\partial_{i'}+\sum_{q,\,r\in Y_1\setminus
i'}c_{qr}x_q\partial_r.\] For $ j'\in Y_1\setminus i',$ it is easy to see that
$[\mathrm{T_H}(x^{(a\varepsilon_
i)}x_{i'}),\,\mathrm{T_H}(x_ix_{j'})]=-\mathrm{T_H}(x^{(a\varepsilon_i)}x_{j'}).$
Applying $\phi$ we obtain that
\[-c_{i'i'}x_{j'}\partial_{i'}+\sum_{q\in Y_1\setminus i'}c_{qj'}x_q\partial_{i'}=0.\]
It follows that
$c_{i'i'}=c_{j'j'}$ and  $ c_{qj'}=0$ for $ q\in Y_1\setminus
\{i',j'\}$ and $ j'\in Y_1\setminus i'.$
Then
$$\phi(\mathrm{T_H}(x^{(a\varepsilon_i)}x_{i'}))=\sum_{r\in
Y_1\setminus\{i',j'\}}c_{j'r}x_{j'}\partial_r+c_{i'i'}x_{i'}\partial_{i'}
+c_{i'i'}x_{j'}\partial_{j'}.$$
For  $l\in Y_1\setminus\{i',\,j'\},$ we have
$$[\mathrm{T_H}(x^{(a\varepsilon_i)}x_{i'}),\,\mathrm{T_H}(x_ix_l)]=-\mathrm{T_H}(x^{(a\varepsilon_i)}x_l).$$
Applying $\phi$ we obtain that
  $c_{j'l}x_{j'}\partial_{i'}-c_{i'i'}x_l\partial_{i'}=0.$
This implies that $c_{i'i'}=0$ and $c_{j'l}=0$ for all
$l\in Y_1\setminus\{i',j'\}.$
Therefore  $\phi(\mathrm{T_H}(x^{(a\varepsilon_i)}x_{i'}))=0.$\\

\noindent Case (ii):
$a\equiv 2\pmod{p}.$ Then
 \[\phi(\mathrm{T_H}(x^{(a\varepsilon_i)}x_{i'}))=\sum_{q\in
Y_1\setminus i'}c_{qi'}x_q\partial_{i'}.\]
Note that
$[\mathrm{T_H}(x^{(a\varepsilon_i)}x_{i'}),\,\Delta_{j}]=0 $
for $
j\in Y_0\setminus i.$ Similarly, applying $\phi$ one gets  $c_{j'i'}=0$ for $ j\in
Y_0\setminus i $ and therefore,
 $\phi(\mathrm{T_H}(x^{(a\varepsilon_i)}x_{i'}))=0.$\\

\noindent Case (iii):
$a\not\equiv 1, 2\pmod{p}.$
Clearly, $\phi(\mathrm{T_H}(x^{(a\varepsilon_i)}x_{i'}))=0.$
\end{proof}
\begin{proposition}\label{ht4.6}
Suppose  $t\in \mathbb{N}$ is not any $p$-power.
Then $\mathrm{Der}_{-t}(\mathcal{HO},\mathcal{W})=0.$
\end{proposition}
\begin{proof} Let $\phi\in
\mathrm{Der}_{-t}(\mathcal{HO},\mathcal{W}).$
We first consider the case $t\not\equiv 0\pmod{p}.$
Direct computation shows
$$[\Delta_{i},\mathrm{T_H}(x^{(t\varepsilon_i)}x_k)]=
(\delta_{k,i'}-t)\mathrm{T_H}(x^{(t\varepsilon_i)}x_k).$$
Applying $\phi$ one gets
$$[\Delta_{i},\,\phi(\mathrm{T_H}(x^{(t\varepsilon_i)}x_k))]=
(\delta_{k,i'}-t)\phi(\mathrm{T_H}(x^{(t\varepsilon_i)}x_k)).$$
Since $\phi(\mathrm{T_H}(x^{(t\varepsilon_i)}x_k))\in
\mathcal{HO}_{-1},$
one may assume that
\[\phi(\mathrm{T_H}(x^{(t\varepsilon_i)}x_k))=\sum_{r\in
Y_0}a_r\partial_r\quad\mbox{where}\,a_r\in \mathbb{F}.\]
Then $$\Big[\Delta_{i},\,\sum_{r\in
Y_0}a_r\partial_r\Big]=a_i\partial_i=(\delta_{k,i'}-t)\sum_{r\in Y_0}a_r\partial_r.$$
If $k\neq i',$ then $(t+1)a_i=0 $ and $ta_r=0$ for $ r\in Y_0\setminus i.$
Since $t\not\equiv 0\pmod{p},$  we have $a_r=0$ for $ r\in Y_0\setminus i.$
If $t\not\equiv -1\pmod{p},$ then $a_i=0 $ and therefore,
 $\phi(\mathrm{T_H}(x^{(t\varepsilon_i)}x_k))=0.$
 If $t\equiv
-1\pmod{p}$, then $\phi(\mathrm{T_H}(x^{(t\varepsilon_i)}x_k))=a_i\partial_i.$
Note that
$$[\mathrm{T_H}(x^{(t\varepsilon_i)}x_k),\Delta_{k'}]
=-\mathrm{T_H}(x^{(t\varepsilon_i)}x_k).$$
Applying $\phi$ one gets  $[a_i\partial_i,\,x_k\partial_k-x_{k'}\partial_{k'}]=-a_i\partial_i $  and therefore, $
a_i=0.$
Hence  $$\phi(\mathrm{T_H}(x^{(t\varepsilon_i)}x_k))=0.$$
If $k=i',$ then $ta_i=0 $ and $(t-1)a_r=0$ for $ r\in Y_0\setminus i.$
Since $t\not\equiv 0\pmod{p},$ we have $a_i=0.$ If $t\not\equiv 1\pmod{p},$
then $a_r=0$ for $ r\in Y_0\setminus i $ and hence
 $$\phi(\mathrm{T_H}(x^{(t\varepsilon_i)}x_{i'}))=0.$$
If $t\equiv 1\pmod{p},$
then  $\phi(\mathrm{T_H}(x^{(t\varepsilon_i)}x_{i'}))=\sum_{r\in
Y_0\setminus i}a_r\partial_r.$ Note that
 $[\mathrm{T_H}(x^{(t\varepsilon_i)}x_{i'}),\,\Delta_{j}]=0 $
 for $j\in Y_0\setminus i.$
Applying $\phi$ one gets
   $-a_j\partial_j=0 $ and $a_j=0.$
Thus
 $\phi(\mathrm{T_H}(x^{(t\varepsilon_i)}x_{i'}))=0.$

 Let us consider the case $t\equiv 0\pmod{p}.$
 Write $t=\sum_{s=1}^ra_sp^s $ where $ 0\leq
a_s<p$ and $a_r\neq 0.$ Since
$$\mathrm{zd}(\mathrm{T_H}(x^{((t-p^r+1)\varepsilon_i)}x_k))=t-p^r<t-2,$$
$$\mathrm{zd}(\mathrm{T_H}(x^{(p^r\varepsilon_i)}x_{i'}))=p^r-1<t-2,$$
we have
$$\phi(\mathrm{T_H}(x^{((t-p^r+1)\varepsilon_i)}x_k))=\phi(\mathrm{T_H}(x^{(p^r\varepsilon_i)}x_{i'}))=0.$$
Then
\begin{equation*}
[\mathrm{T_H}(x^{((t-p^r+1)\varepsilon_i)}x_k),\,
\mathrm{T_H}(x^{(p^r\varepsilon_i)}x_{i'})]
=\Big[{t\choose {p^r}}-\delta_{k,i'}{t\choose{p^r-1}}\Big
]\mathrm{T_H}(x^{(t\varepsilon_i)}x_k).
\end{equation*}
Note that ${t\choose{p^r-1}}\equiv 0\pmod{p} $ and ${t\choose{p^r}}\not\equiv
0\pmod{p}.$ Applying $\phi$ we have
$$\phi(\mathrm{T_H}(x^{(t\varepsilon_i)}x_k))=0 \quad\mbox{for all}\, i\in
Y_0,\,k\in Y_1.$$ By Lemma \ref{ht4.5}, $\phi=0.$ This proves
 $\mathrm{Der}_{-t}(\mathcal{HO},\mathcal{W})=0.$
\end{proof}
\begin{proposition}\label{ht4.7} Let $t=p^r$ for some $r\in \mathbb{N}.$
Then
$\mathrm{Der}_{-t}(\mathcal{HO},\mathcal{W})
=\mathrm{span}_\mathbb{F}\{(\mathrm{ad}\partial_i)^t|i\in
Y_0\}.$
\end{proposition}
\begin{proof}
Let $\phi\in \mathrm{Der}_{-t}(\mathcal{HO},\mathcal{W}).$
Then $\phi(\mathrm{T_H}(x^{(t\varepsilon_i)}x_k))\in
\mathcal{HO}_{-1} $ for $i\in Y_0$ and $k\in Y_1.$
Assume that
$$\phi(\mathrm{T_H}(x^{(t\varepsilon_i)}x_k))=\sum_{r\in
Y_0}a_{ikr}\partial_r\quad\mbox{where}\,\;  a_{ikr}\in \mathbb{F}.
$$
Compute
$$[\mathrm{T_H}(x^{(t\varepsilon_i)}x_k),\Delta_{i}]
=(t-\delta_{k,i'})\mathrm{T_H}(x^{(t\varepsilon_i)}x_k).$$
Applying $\phi$ one gets
\[-a_{iki}\partial_i=(t-\delta_{k,i'})\sum_{r\in
Y_0}a_{ikr}\partial_r=-\delta_{ki'}\sum_{r\in Y_0}a_{ikr}\partial_r.\]
If $i'=k,$ then $a_{ikr}=0$ for $ r\in Y_0\setminus i $ and hence
 $$\phi(\mathrm{T_H}(x^{(t\varepsilon_i)}x_{i'}))=a_{ii'i}\partial_i.$$
If  $i'\neq k,$ then $a_{iki}=0 $ and hence
 \[\phi(\mathrm{T_H}(x^{(t\varepsilon_i)}x_k))=\sum_{r\in
Y_0\setminus i}a_{ikr}\partial_r.\]
Note that  $[\mathrm{T_H}(x^{(t\varepsilon_i)}x_k),\Delta_{k'})]
=-\mathrm{T_H}(x^{(t\varepsilon_i)}x_k).$
Applying $\phi,$  one gets
$a_{ikr}=0$ for $r\in Y_0\setminus k' $ and therefore,
 $$\phi(\mathrm{T_H}(x^{(t\varepsilon_i)}x_k))=a_{ikk'}\partial_{k'}.$$
Note that
$$[\mathrm{T_H}(x^{(t\varepsilon_i)}x_{i'}),\,\mathrm{T_H}(x_ix_k)]
=-\mathrm{T_H}(x^{(t\varepsilon_i)}x_k).$$
Applying $\phi$, we have
 $[a_{ii'i}\partial_i,x_k\partial_{i'}-x_i\partial_{k'}]=-a_{ikk'}\partial_{k'} $  and it follows that
   $a_{ii'i}=a_{ikk'}$ for $
k\neq i'.$ Put $a_i:=a_{ii'i}=a_{ikk'} $
and \[\psi:=\phi+\sum_{r\in Y_0}a_r(\mathrm{ad}\partial_r)^t.\]
Then $\psi\in
\mathrm{Der}_{-t}(\mathcal{HO},\mathcal{W}) $ and one may easily verify that
$\psi(\mathrm{T_H}(x^{(t\varepsilon_i)}x_k))
=0.$
By Lemma \ref{ht4.5}, $\psi=0 $ and
 $\phi=-\sum_{r\in Y_0}a_r(\mathrm{ad}\partial_r)^t\in
\mathrm{span}_\mathbb{F}\{(\mathrm{ad}\partial_i)^t|i \in Y_0\}.$
\end{proof}

Assembling the main results in Section 2 and 3, we are able to describe explicitly the derivation
space from $\mathcal{HO}$ into $\mathcal{W}:$

\begin{theorem}\label{ht4.8}
  $\mathrm{Der}(\mathcal{HO},\mathcal{W})
  =\mathrm{ad}\mathcal{W}+\mathrm{span}_{\mathbb{F}}\{(\mathrm{ad}\partial_{i})^{p^{k_{i}}}\mid i\in
Y_0,\,1\leq k_i<t_i\}.$
\end{theorem}
\begin{proof}
  This is a direct consequence of Propositions \ref{ht3.8}, \ref{ht4.4}, \ref{ht4.6} and \ref{ht4.7}.
\end{proof}

\section{The derivation algebra of $\mathcal{HO}$}
In this section, using the results obtained in Sections 2 and 3, we shall determine the derivation algebra and
outerderivation algebra of for the even part of the Hamiltonian superalgebra.  Note that
by \cite[Proposition 20, p. 202]{lzw},
$$\mathcal{HO}=\overline{\mathcal{HO}}=\Big\{\sum_{i=1}^{2n}a_i\partial_i\in
 \mathcal{W}\,\Big|\,\partial_i(a_{j'})=(-1)^{\mu(i)\mu(j)+\mu(i)+\mu(j)}\partial_j(a_{i'}),
 \, \,i,j\in
 Y\Big\}.$$

\begin{proposition}\label{ht5.1} Let $t$ be a positive integer. Then
 $\mathrm{Der}_t(\mathcal{HO})=\mathrm{ad}(\mathcal{HO}_t).$
\end{proposition}
 \begin{proof} Let $\phi\in \mathrm{Der}_t(\mathcal{HO}).$ Then
 $\phi\in
 \mathrm{Der}_t(\mathcal{HO},\mathcal{W}).$
 By Proposition \ref{ht3.8}, there is $D\in
 \mathcal{W}_t $ such that
 $\phi=\mathrm{ad}D\in \mathrm{Der}_t(\mathcal{HO}).$ Clearly,  $D\in
 \mathrm{Nor}_{\mathcal{W}}(\mathcal{HO})_t.$
 Let $D=\sum_{r\in Y}g_r\partial_r, $ where $g_r\in
 \mathcal{O}(n,n;\underline{t})_{t+1}.$ Then
\begin{equation}\label{he5.1}
[\partial_i,\,D]=\Big[\partial_i,\,\sum_{r\in Y}g_r\partial_r\Big]=\sum_{r\in
Y}\partial_i(g_r)\partial_r\in \mathcal{HO}\quad\mbox{for}\,i\in Y_0.
\end{equation}
Fix $i\in Y_0 $ and put  $a_k:=\partial_i(g_k) $ for $k\in
Y.$
By the definition of $\mathcal{HO}$,  $\partial_i(a_i)=-\partial_{i'}(a_{i'}),$
 that is, $\partial_i(\partial_i(g_i))=-\partial_{i'}(\partial_i(g_{i'})).$ Therefore,
 \begin{equation}
 \partial_i(\partial_i(g_i)+\partial_{i'}(g_{i'}))=0.\label{he5.2}
 \end{equation}
Similarly, we have $\partial_{i'}(a_i)=-\partial_{i'}(a_i), $ that is,  $\partial_{i'}(a_i)=0.$ Hence
\begin{equation}
\partial_{i'}(\partial_i(g_i))=0.\label{he5.3}
 \end{equation}
 For any $r\in
Y\setminus\{i,i'\},$ $\partial_i(a_r)=(-1)^{\mu(r')}\partial_{r'}(a_{i'}),$
that is,  $\partial_i(\partial_i(g_r))=(-1)^{\mu(r')}\partial_{r'}(\partial_i(g_{i'})).$
It follows that
\begin{equation}
\partial_i(\partial_i(g_r)-(-1)^{\mu(r')}\partial_{r'}(g_{i'}))=0.\label{he5.4}
 \end{equation}
For $r\in Y\setminus\{i,i'\},$ we have $\partial_{i'}(a_r)=-\partial_{r'}(a_i),$
that is, $\partial_{i'}(\partial_i(g_r))=-\partial_{r'}(\partial_i(g_i)).$  It follows that
\begin{equation} \partial_i(\partial_{i'}(g_r)+\partial_{r'}(g_i))=0.\label{he5.5}\end{equation}
 For $r\in
Y\setminus\{i,i'\},$ we have $\partial_r(a_r)=-\partial_{r'}(a_{r'}),$
that is, $\partial_r(\partial_i(g_r))=-\partial_{r'}(\partial_i(g_{r'})).$
It follows that
\begin{equation}
 \partial_i(\partial_r(g_r)+\partial_{r'}(g_{r'}))=0.\label{he5.6}
 \end{equation}
Note that
\begin{equation}
[\Delta_{i},\,D]
=( \Delta_{i}(g_{i'})-g_{i'})\partial_{i'}+(\Delta_{i}(g_{i})+g_i)\partial_i
  +\sum_{r\in
 Y\setminus
 \{i,i'\}}\Delta_{i}(g_r)\partial_r\in \mathcal{HO}.
 \label{he5.7}\end{equation}
 For $k\in Y,$ denote by $b_k$ the coefficient  of $\partial_k$.
 Then by the definition of $\mathcal{HO},$ we have $\partial_{i'}(b_i)=-\partial_{i'}(b_i),$
 and therefore,
$$\partial_{i'}(\Delta_{i}(g_i)+g_i)=0.$$
It follows that
$$\partial_{i'}(g_i)-x_i\partial_{i'}\partial_i(g_i)+\partial_{i'}(g_i)=0.$$
By  (\ref{he5.3}), $\partial_{i'}(g_i)+\partial_{i'}(g_i)=0.$
Similarly, for  $r\in Y\setminus\{i,i'\},$ we have $\partial_{i'}(b_r)=-\partial_{r'}(b_i),$
that is,
$$\partial_{i'}(\Delta_{i}(g_r))=-\partial_{r'}(\Delta_{i}(g_i)+g_i).$$
Computing from the equation above one gets
$$\partial_{i'}(g_r)+\partial_{r'}(g_i)=-\Delta_{i}(\partial_{i'}(g_r)+\partial_{r'}(g_i)).$$
 By (\ref{he5.5}), we have
 $x_{i'}\partial_{i'}(\partial_{i'}(g_r)+\partial_{r'}(g_i))=-(\partial_{i'}(g_r)+\partial_{r'}(g_i)). $
It follows that
\begin{equation}
\partial_{i'}(g_r)+\partial_{r'}(g_i)=0.\label{he5.8}
\end{equation}
For $r\in Y\setminus\{i,i'\},\,\partial_r(b_r)=-\partial_{r'}(b_{r'}), $ and then
$$\partial_r (\Delta_{i}(g_r)) =-\partial_{r'} (\Delta_{i}(g_{r'})).$$
Computing from the equation above one gets
$$x_{i'}\partial_{i'}(\partial_r(g_r)+\partial_{r'}(g_{r'}))=x_i\partial_i(\partial_r(g_r)+\partial_{r'}(g_{r'})).$$
By (\ref{he5.6}),
 $x_{i'}\partial_{i'}(\partial_r(g_r)+\partial_{r'}(g_{r'}))=0.$
Consequently,
\begin{equation}
\partial_{i'}(\partial_r(g_r)+\partial_{r'}(g_{r'}))=0.\label{he5.9}
\end{equation}
For $r\in
Y\setminus\{i,i'\},$ we have $\partial_i(b_r)=(-1)^{\mu(r')}\partial_{r'}(b_{i'}),$
that is,
$$\partial_i (\Delta_{i}(g_r))
=(-1)^{\mu(r')}\partial_{r'}(\Delta_{i}(g_{i'})-g_{i'}).$$
Similarly, one computes
$$ \Delta_{i} (\partial_i(g_r)-(-1)^{\mu(r')}\partial_{r'}(g_{i'}))
=\partial_i(g_r)-(-1)^{\mu(r')}\partial_{r'}(g_{i'}).$$
By (\ref{he5.4}),
\begin{equation}
x_{i'}\partial_{i'}(\partial_i(g_r)-(-1)^{\mu(r')}\partial_{r'}(g_{i'}))
=\partial_i(g_r)-(-1)^{\mu(r')}\partial_{r'}(g_{i'}).\label{he5.10}
\end{equation}
Similar to the proof of \cite[Lemma 6]{lzw}, one may show that
\begin{equation} \partial_i(g_{r'})=\partial_r(g_{i'}) \quad\mbox{for}\,\,i,r\in
Y_0 \,\;\mbox{with}\,\; i\not=r.\label{he5.11}
\end{equation}

We next show that $\partial_{r'}(g_r)=-\partial_{r'}(g_r) $ for $r\in Y_0\setminus i.$
By (\ref{he5.1}), $\partial_{r'}(a_r)=-\partial_{r'}(a_r) $ and then
 $\partial_{r'}(\partial_i(g_r))=0.$ Hence $\partial_i(\partial_{r'}(g_r))=0.$
 By (\ref{he5.7}), $\partial_{r'}(b_r)=-\partial_{r'}(b_r) $ and then
  $\partial_{r'}\Delta_{i}(g_r)=0.$
It follows that
$x_{i'}\partial_{i'}\partial_{r'}(g_r)=0  $ and then $\partial_{i'}\partial_{r'}(g_r)=0,$
since we have shown that $\partial_i(\partial_{r'}(g_r))=0.$
Clearly, $\partial_{r'}(\partial_{r'}(g_r))=0$. Since $D\in
\mathrm{Nor}_\mathcal{W}(\mathcal{HO})_t, $ there is
$f_r\in
\mathcal{O}(n,n;\underline{t})_{\overline{1}} $
such that $[\partial_r,\,D]=\mathrm{T_H}(f_r)$ for $  r\in Y_0\setminus i.$
Consequently,
$$\sum_{j\in Y}\partial_r(g_j)\partial_j=\sum_{j\in
Y}(-1)^{\mu(j)}\partial_j(f_r)\partial_{j'}.$$ Therefore,
$\partial_r(g_r)=-\partial_{r'}(f_r)$ for $ r\in
Y_0\setminus i $ and then
 $\partial_r\partial_{r'}(g_r)=\partial_{r'}\partial_r(g_r)=-\partial_{r'}\partial_{r'}(f_r)=0.$ Thus we have
\begin{equation}
\partial_{r'}(g_r)=0.\label{he5.12}
\end{equation}
It follows that $\partial_{r'}(g_r)=-\partial_{r'}(g_r).$
We have to show that for $t>1$,
\begin{equation} \partial_i(g_r)=-\partial_{r'}(g_{i'})\quad\mbox{for}\,\; r\in Y_0\setminus
i.\label{he5.13}
\end{equation}
Form  (\ref{he5.4}) one can see  $\partial_i(\partial_i(g_r)+\partial_{r'}(g_{i'}))=0.$
Then by  (\ref{he5.12}),
 $$\partial_{r'}(\partial_i(g_r)+\partial_{r'}(g_{i'}))
 =\partial_i(\partial_{r'}(g_r))+\partial_{r'}\partial_{r'}(g_{i'})=0\quad \mbox{for all}\; r\in Y_0\setminus i.$$
By (\ref{he5.11}) and (\ref{he5.6}),  for $r\in Y_0\setminus i,$ one may compute
 \begin{eqnarray*}
 \partial_r(\partial_i(g_r)+\partial_{r'}(g_{i'}))
 =0.
\end{eqnarray*}
Note also that, by (\ref{he5.10}),
$x_{i'}\partial_{i'}(\partial_i(g_r)+\partial_{r'}(g_{i'}))=\partial_i(g_r)+\partial_{r'}(g_{i'}).$ Then one may see that
$\partial_i(g_r)+\partial_{r'}(g_{i'})$ is of the form $\lambda x_{i'} $ where
$\lambda\in \mathbb{F}.$ If $t>1,$
then  $\partial_i(g_r)+\partial_{r'}(g_{i'})=0.$ In the following we show that
  $\partial_r(g_r)=-\partial_{r'}(g_{r'})$ for $t>1.$
By (\ref{he5.6}), $\partial_i(\partial_r(g_r)+\partial_{r'}(g_{r'}))=0.$
By (\ref{he5.9}), $\partial_{i'}(\partial_r(g_r)+\partial_{r'}(g_{r'}))=0.$ By  (\ref{he5.12}),
$$\partial_{r'}(\partial_r(g_r)+\partial_{r'}(g_{r'}))=\partial_r\partial_{r'}(g_r)+\partial_{r'}\partial_{r'}(g_{r'})=0.$$
Since $r\in Y_0\setminus i,$
 $\partial_r(g_r)+\partial_{r'}(g_{r'})$ is of the form $\lambda x_r $
 where $\lambda\in \mathbb{F}.$ If $t>1,$
then $$\partial_r(g_r)+\partial_{r'}(g_{r'})=0.$$
Let us show that
$\partial_i(g_r)=-\partial_{r'}(g_{i'})$ and $ \partial_r(g_r)=-\partial_{r'}(g_{r'}) $ for $t=1$.
Let $D:=\sum_{j\in Y}g_j\partial_j,$ where $g_j=\sum_{k,l\in
Y}c_{jkl}x_kx_l.$ Then for arbitrary $i\in Y_0,$
\begin{eqnarray*}[\partial_i,\,D]&=&\Big[\partial_i,\,\sum_{k,l,j\in
Y}c_{jkl}x_kx_l\partial_j\Big]\\&=&\sum_{l,j\in
Y}c_{jil}x_l\partial_j+\sum_{k,j\in Y}c_{jki}x_k\partial_j\\&=&\sum_{l,j\in
Y}(c_{jil}+c_{jli})x_l\partial_j\\&=&\sum_{j\in Y}\Big(\sum_{l\in
Y}(c_{jil}+c_{jli})x_l\Big)\partial_j\in \mathcal{HO}.
\end{eqnarray*}
For fixed  $i\in Y_{0}, $  put  $d_j:= \sum_{l\in
Y}(c_{jil}+c_{jli})x_l $ for $j\in Y.$ Since
$\partial_r(d_r)=-\partial_{r'}(d_{r'}) $  for $r\in Y_0\setminus i,$ one gets
\begin{equation}
c_{rir}+c_{rri}=-(c_{r'ir'}+c_{r'r'i}).\label{he5.14}
\end{equation}
For $r\in Y_0\setminus i,$ we have $ \partial_i(d_r)=-\partial_{r'}(d_{i'}).$
It follows that
\begin{equation} 2c_{rii}=-(c_{i'ir'}+c_{i'r'i}).\label{he5.15}
\end{equation}
For arbitrary $r\in Y_0\setminus i,$ similarly, we have
\begin{eqnarray*}
[\partial_r,\,D]=\sum_{j\in Y}\Big(\sum_{l\in
Y}(c_{jrl}+c_{jlr})x_l\Big)\partial_j\in \mathcal{HO}.
\end{eqnarray*}
Put  $e_j:=\sum_{l\in
Y}(c_{jrl}+c_{jlr})x_l $ for $j\in Y.$ For $r\in Y_0\setminus i, $ we have $
\partial_r(e_r)=-\partial_{r'}(e_{r'}).$ It follows that
\begin{equation}
2c_{rrr} =-(c_{r'rr'}+c_{r'r'r}).\label{he5.16}
\end{equation}
For arbitrary $r\in Y_0\setminus i,$ we have $ \partial_i(e_r)=-\partial_{r'}(e_{i'}) $
and it follows that
\begin{equation} c_{rri}+c_{rir}=-(c_{i'rr'}+c_{i'r'r}).\label{he5.17}
\end{equation}
By (\ref{he5.14}) and (\ref{he5.16}), we have
\begin{eqnarray*}\partial_r(g_r)&=&\partial_r\Big(\sum_{k,l\in
Y}c_{rkl}x_kx_l\Big)\\&=&\sum_{l\in Y_0}c_{rrl}x_l+\sum_{k\in
Y_0}c_{rkr}x_k\\&=&\sum_{l\in
Y_0}(c_{rrl}+c_{rlr})x_l\\&=&-\sum_{l\in
Y_0}(c_{r'r'l}+c_{r'lr'})x_l\\&=&-\sum_{l\in
Y_0}c_{r'r'l}x_l-\sum_{k\in
Y_0}c_{r'kr'}x_k\\&=&-\partial_{r'}\Big(\sum_{k,l\in
Y}c_{r'kl}x_kx_l\Big)\\&=&-\partial_{r'}(g_{r'}).
\end{eqnarray*}
 Similarly, by (\ref{he5.15}) and (\ref{he5.17}), one may verify that
\begin{equation}
\partial_i(g_r)=-\partial_{r'}(g_{i'}).\label{he5.18}
\end{equation}
By (\ref{he5.3}), $\partial_{i'}(\partial_i(g_i))=0 $
 and then
  $\partial_{i'}(\partial_i(g_i)+\partial_{i'}(g_{i'}))=0.$  By (\ref{he5.2}),
$$\partial_i(\partial_i(g_i)+\partial_{i'}(g_{i'}))=0.$$
For $r\in Y_0\setminus i,$ by (\ref{he5.8}) and (\ref{he5.11}), we have
$\partial_r(\partial_i(g_i)+\partial_{i'}(g_{i'}))
 =0.$
For $r\in Y_0\setminus i,$ by  (\ref{he5.8}), (\ref{he5.13}) and
(\ref{he5.18}), we have
$\partial_{r'}(\partial_i(g_i)+\partial_{i'}(g_{i'}))=0.$ Consequently, $\partial_i(g_i)=-\partial_{i'}(g_{i'}).$
By the definition of $\mathcal{HO},$   $D\in \mathcal{HO}_t $ and $\phi\in
\mathrm{ad}\mathcal{HO}_t.$
This completes the proof.
\end{proof}

Put $\Gamma:=\sum_{i=1}^nx_{i'}\partial_{i'}.$

\begin{proposition}\label{ht5.2}
$\mathrm{Der}_0(\mathcal{HO})=\mathrm{ad}(\mathcal{HO}+\mathbb{F}\;\Gamma)_0.$
\end{proposition}
\begin{proof} It is easily seen that
$\mathrm{ad}(\mathcal{HO}+\mathbb{F}\,\Gamma)_0\subset
 \mathrm{Der}_0{(\mathcal{HO})}.$
 Let $\phi\in \mathrm{Der}_0{(\mathcal{HO})}.$ Then
 $\phi\in
 \mathrm{Der}_0{(\mathcal{HO},\mathcal{W})}.$ By Proposition \ref{ht3.8}, there is
 $D\in
 \mathcal{W}_0 $ such that
 $\phi=\mathrm{ad}D\in \mathrm{Der}_0(\mathcal{HO}) .$ Clearly, $D\in
 \mathrm{Nor}_{\mathcal{W}}(\mathcal{HO})_0.$
Let $D:=\sum_{k,l\in Y}c_{kl}x_k\partial_l.$ For  $i\in Y_0,j\in
 Y_1,$ direct computation shows that
 \begin{eqnarray*}
 [\mathrm{T_H}(x_ix_j),\,D]&=&\Big[x_j\partial_{i'}-x_i\partial_{j'},\,\sum_{k,l\in
 Y}c_{kl}x_k\partial_l\Big]\\
 &=&\Big(c_{i'i'}x_j-\sum_{k\in Y_1}c_{kj}x_k\Big)\partial_{i'}+\Big(\sum_{k\in
 Y_0}c_{ki}x_k-c_{j'j'}x_i\Big)\partial_{j'}\\&+&\sum_{l\in Y_1\setminus
 {i'}}c_{i'l}x_j\partial_l-\sum_{r\in Y_0\setminus
 {j'}}c_{j'r}x_i\partial_r.\end{eqnarray*}
 Note that $[\mathrm{T_H}(x_ix_j),D]\in \overline{\mathcal{HO}}=\mathcal{HO}.$
 Denote by $a_k $ the coefficient of $\partial_k$ in the equation above. Then
 $\partial_i(a_{j'})=-\partial_j(a_{i'}).$  It follows that
 $c_{ii}+c_{i'i'}=c_{jj}+c_{j'j'} $ for all $i \in Y_0 $ and $j\in Y_1.$
 Let $c:=c_{ii}+c_{i'i'}$
 for $i\in Y_0.$
 In the same way one may obtain that
 $\partial_j(a_l)=-\partial_{l'}(a_{j'})$ for $ l\in Y_1\setminus i'.$
 It follows that
 $c_{i'l}=-c_{l'i} $  for $ l\in Y_1\setminus i'.$
 Since $D\in \mathcal{W}_0,$  $c_{ij}=0$ for $ i\in Y_0,j\in Y_1 $ and then
 $D=\sum_{k,l\in Y}c_{kl}x_k\partial_l,$ where $\mu(k)= \mu(l).$
Then
 \begin{eqnarray*}D&=&\sum_{k,l\in
 Y_1}c_{kl}x_k\partial_l+\sum_{i,j\in Y_0}c_{ij}x_i\partial_j\\
 &=&\sum_{k\in Y_1}c_{kk}x_k\partial_k+\sum_{i\in
 Y_0}c_{ii}x_i\partial_i+\sum_{k,l\in Y_1\atop k\neq l}c_{kl}x_k\partial_l+\sum_{i,j\in
 Y_0\atop i\neq j}c_{ij}x_i\partial_j\\
 &=&\sum_{i=1}^n(c_{ii}x_i\partial_i+c_{i'i'}x_{i'}\partial_{i'})+\sum_{k,l=1\atop k\neq
 l}^n(c_{kl}x_k\partial_l+c_{l'k'}x_{l'}\partial_{k'})\\
 &=&\sum_{i=1}^n(c_{ii}x_i\partial_i-c_{ii}x_{i'}\partial_{i'}+c_{ii}x_{i'}\partial_{i'}+c_{i'i'}x_{i'}\partial_{i'})+\sum_{k,l=1\atop k\neq
 l}^n(c_{kl}x_k\partial_l-c_{kl}x_{l'}\partial_{k'})\\
 &=&-\sum_{i=1}^nc_{ii}\Delta_{i}+c\Gamma-\sum_{k,l=1\atop k\neq
 l}^nc_{kl}\mathrm{T_H}(x_kx_{l'})\in
 \mathcal{HO}+\mathbb{F}\,\Gamma.\end{eqnarray*}
 Hence,
 $\phi=\mathrm{ad}D\in \mathrm{ad}(\mathcal{HO}+\mathbb{F}\,\Gamma)_0.$
 \end{proof}
\begin{proposition}\label{ht5.3}
 $\mathrm{Der}_{-1}(\mathcal{HO})=\mathrm{ad}\mathcal{HO}_{-1}.$
\end{proposition}
\begin{proof} This is a direct consequence of  Proposition \ref{ht4.4}.
\end{proof}
  \begin{proposition}\label{ht5.4} If $m>1 $ is not any $p$-power. Then
$\mathrm{Der}_{-m}(\mathcal{HO})=0.$ If $m=p^r$ for some positive integer $r,$ then
$\mathrm{Der}_{-m}(\mathcal{HO})=\mathrm{span}_\mathbb{F}\{(\mathrm{ad}\partial_i)^t|i\in
Y_0\}.$
\end{proposition}
\begin{proof} This is a direct consequence of  Propositions \ref{ht4.6} and \ref{ht4.7}.
\end{proof}

Recall the notation $\underline{t}=(t_{1},\cdots,t_{n})$ stands for a fixed $n$-tuple of
positive integers defining the finite-dimensional odd Hamiltonian superalgebra (see Section 1).
\begin{theorem}\label{ht5.5}
$\mathrm{Der}(\mathcal{HO})=\mathrm{ad}(\mathcal{HO}+\mathbb{F}\, \Gamma)\bigoplus
\mathrm{span}_\mathbb{F}\{(\mathrm{ad}\partial_i)^{p^{k_i}}|i\in
Y_0,\,1\leq k_i<t_i\}.$
\end{theorem}
\begin{proof} This is a direct consequence of  Propositions \ref{ht5.1}--  \ref{ht5.4}.
\end{proof}
Put $\mathcal{J}:=\mathrm{span}_\mathbb{F}\{\partial_i^{p^{r_i}}|i\in
Y_0,1\leq r_i<t_i\}.$
Then
$\mathcal{J}$ is an abelian sub-Lie algebra of
$\mathrm{Der}(\mathcal{O}(n,n;\underline{t}))$ of dimension $(\sum_{i\in
Y_0}t_i-n)$, and $[\mathcal{J},\Gamma]=0. $
Furthermore, by Theorem  \ref{ht5.5}  we have

\begin{theorem}\label{ht5.7}
The outer derivation algebra of  $\mathcal{HO}$  is an abelian Lie algebra of dimension
$\sum_{i\in
Y_0}t_i-n+1.$
\end{theorem}

We conclude this section by the following corollary.
\begin{corollary}  \label{ht5.8}
$\mathrm{dim}_\mathbb{F}(\mathrm{Der}(\mathcal{HO}))=2^{n-1}p^{\sum_{i\in
Y_0}t_i}+\sum_{i\in Y_0}t_i-n.$
\end{corollary}

\begin{proof} By  \cite[Theorem 2.5]{lz2},
$\mathrm{dim}_\mathbb{F}HO=2^np^{\sum_{i\in Y_0}t_i}-1.$
Then $\mathrm{dim}_\mathbb{F}\mathcal{HO}=2^{n-1}p^{\sum_{i\in
Y_0}t_i}-1.$
Let $a\in C(\mathcal{HO}).$ Since
$C_\mathcal{HO}(\mathcal{HO}_{-1})\subset\mathcal{G},$  we have
$a\in
\mathcal{G}\bigcap\mathcal{HO}.$ Let  $a=\sum_{r\in Y_0,u\in
\mathbb{B}^0}c_{ur}x^u\partial_r, $ where $c_{ur}\in \mathbb{F}.$ Then
 $\Big[\Delta_{i},\,\sum_{r\in Y_0,u\in
\mathbb{B}^0}c_{ur}x^u\partial_r\Big]=0.$ It follows that
$$\delta_{i'\in u}\sum_{r\in
Y_0,u\in \mathbb{B}^0}c_{ur}x^u\partial_r+\sum_{u\in
\mathbb{B}^0}c_{ui}x^u\partial_i=0.$$ For any fixed $u,$ if $i'\in u,$
then $c_{ur}=0$ for $ r\in Y_0;$ if $i'\not\in u,$ then $c_{ui}=0.$
Hence, $a=0,$ proving $C(\mathcal{HO})=0.$
Now our conclusion follows from Theorem \ref{ht5.7}.
\end{proof}

\vspace{0.5cm}

\end{document}